\documentclass[fleqn]{article}
\usepackage{graphicx}
\usepackage{amsmath}
\usepackage{amsthm}
\usepackage{hyperref}
\usepackage{amssymb}
\usepackage{enumitem}
\usepackage{stmaryrd}
\usepackage{mathtools}
\usepackage[nonewpage]{imakeidx}
\usepackage{relsize}
\usepackage{scalerel}
\usepackage{tikz,tikz-cd}
\usetikzlibrary{arrows,calc,shapes.geometric,positioning}
\usepackage{adjustbox}
\usepackage[normalem]{ulem}
\usepackage{xcolor}

\title{Finite local principal ideal rings}
\author{Matth\'e van der Lee}

\DeclareMathOperator{\Aut}{Aut} 
\DeclareMathOperator{\Gal}{Gal} 
\DeclareMathOperator{\mo}{o} 
\DeclareMathOperator{\ord}{ord} 
\DeclareMathOperator{\res}{res} 

\newcommand{\xb}[1][.75]{\hspace{-#1pt}} 
\newcommand{\xf}[1][1]{\hspace{#1pt}} 

\newcommand\xceil[1]{\ensuremath\ulcorner\xb#1\xb\urcorner} 
\newcommand\xchar{\ensuremath\mathrm{char}} 
\newcommand\xcol{\ensuremath{{\xf[2.5]:\xf[2.5]}}} 
\newcommand{\xexp}[2][6]{\ensuremath\mathrel{\raisebox{#1pt}{$\xb[3.73]{\scriptscriptstyle #2}\xf$}}} 
\newcommand{\xh}[1]{\boldsymbol{#1}} 
\newcommand{\xpe}[1]{\ensuremath p^{\xf#1}\xb} 
\newcommand{\xpef}[1]{\ensuremath\xpe{#1\xb[.5]f}}

\newcommand{\xF}{\ensuremath\mathbb{F}}
\newcommand{\xN}{\ensuremath\mathbb{N}}
\newcommand{\xO}{\ensuremath\mathcal{O}}
\newcommand{\xp}{\ensuremath\mathfrak{p}}
\newcommand{\xP}{\ensuremath\xf\mathfrak{P}}

\newcommand{\xQ}{\ensuremath\mathbb{Q}}
\newcommand{\xZ}{\ensuremath\xf\mathbb{Z}}

\newcommand{\xo}[1]{\ensuremath\overline{#1}}
\newcommand{\xs}[2][i]{\ensuremath\xo{{\smash{#2}\vphantom{#1}}}} 

\newcommand{\xprod}[2]{\ensuremath\Pi_{\,#1}^{\,#2}\xf[.3]}
\newcommand{\xsum}[2]{\ensuremath\Sigma_{\,#1}^{\,#2}\xf[.3]}
\newcommand{\xsize}[1]{\ensuremath|\xf[.3]#1\xf[.3]|}

\newcommand{\xRo}{\ensuremath R_\omega}
\newcommand{\xRb}{\ensuremath R_0[\xf[.5]b\xf[.5]]}
\newcommand{\xu}{\ensuremath k^\times\xb}
\newcommand{\xU}{\ensuremath R^\times\xb} 
\newcommand{\xUo}{\ensuremath \xRo^\times\xb}

\newtheorem{thm}{Theorem}
\newtheorem{prop}{Proposition}
\newtheorem{lem}{Lemma}
\newtheorem{cor}{Corollary}
\newtheorem{appt}{Theorem}

\newtheorem{appp}[appt]{Proposition}
\newtheorem{appl}[appt]{Lemma}

\newtheoremstyle{italicnum}{\topsep}{\topsep}{\upshape}{0pt}{\itshape}{.}{5pt plus 1pt minus 1pt}{\thmname{#1}\thmnumber{\text{ }\itshape#2}\thmnote{(#3)}}

\theoremstyle{italicnum}
\newtheorem{xmp}{Example}

\makeatletter
\tikzset{edge node/.code={\expandafter\def\expandafter\tikz@tonodes \expandafter{\tikz@tonodes #1}}}

\newcommand*{\textlabel}[2]{%
  \edef\@currentlabel{#1}
  \phantomsection
  #1\label{#2}}

\makeatother
\tikzset{subseteq/.style={draw=none,edge node={node [sloped, allow upside down, auto=false]{$\subseteq$}}}} 

\makeindex[columns=1, title=Index of symbols]
\begin{document}
\maketitle

\begin{abstract}
\noindent Every finite local principal ideal ring is the homomorphic image of a discrete valuation ring of a number field, and is determined by five invariants. We present an action of a group, non-commutative in general, on the set of Eisenstein polynomials, of degree matching the ramification index of the ring, over the coefficient ring. The action is defined by taking resultants.
\end{abstract}

\medskip

\begin{tabular}{ll}
\xf[19]\scriptsize{\textbf{Keywords:}} & \scriptsize{Principal ideal rings, Finite commutative rings.}\\
\xf[19]\scriptsize{\textbf{2020 MSC:}} & \scriptsize{13F10, 13M99.}
\end{tabular}
\medskip

\section{Introduction}\label{sec:intro}

Finite local PIRs, principal ideal rings, have been studied by Hungerford (\cite{HU}), who demonstrated, using Cohen's structure theorem for complete local rings, that every finite local PIR is the quotient of a principal ideal domain. The principal ideal domain is itself a quotient of the power series ring over a $v$-ring, for which Hungerford gave explicit relations. Recall that a $v$-\textit{ring}, also known as \textit{Cohen ring}, is a complete DVR (discrete valuation ring) of characteristic 0 with maximal ideal generated by a prime number.\\

Here, Th.\ \ref{thm:main} describes $R$ as the homomorphic image of a DVR of the ring of integers of a number field which has the same residue degree $f$ and ramification index $e$ with respect to the prime number $p$ associated to $R$ as the ring $R$ itself. Th.\ \ref{thm:main2} describes a finite local PIR $R$ in terms of five invariants that determine its isomorphism type. One of the invariants is the set of $p\,$-Eisenstein polynomials of degree $e$ over the coefficient ring of $R$ which have a zero in $R$. This set is closed under an action, defined by taking resultants, of a finite group on the full set of Eisenstein polynomials. The precise formulation is Th.\ \ref{thm:orbits}.\\

The author is indebted to mr. Aurel Page for bringing to his attention that non-isomorphic PIRs exist that cannot be distinguished by the first four invariants, of which version 1 of this paper claimed they already determine the isomorphism type. The present edition is a substantial revision of version 2.  For ease of reference, an appendix on algebraic number theory is included.\\

All rings are commutative with 1. The underlying additive group of a ring $R$ is denoted by $R^{\xf[.5]+}$, and $\xU$ is the group of units. For $x,y\in R$, $x\sim y$ means that $x$ and $y$ are associated ($x\in y\xU$). $\xo{x}$ is the residue class of a ring or group element $x$ modulo an ideal or a subgroup clear from the context. The symbol for an object used throughout the text is boldfaced when the object is (re)defined.

\section{Finite local PIRs}\label{sec:pir}

Let $\xh{R}$\index{R@$R$, a finite local PIR} be a finite local principal ideal ring. $R$ is necessarily Artinian, hence zero-dimensional. Let $\xh{\xp}=\pi R$ be the unique prime ideal\index{P1@$\xh{\xp}$, the prime ideal of $R$}, and $\xh{k}=R/\xp\cong\xF_{p^f}$\xb the residue field\index{K@$k$, the residue field of $R$}, where $\xh{p}$\index{P0@$p$, the residue field characteristic} is a prime number. $\xh{\pi}$\index{YPi@$\pi$, a uniformizer of $R$} is called a \textit{uniformizer} of $R$, and $\xh{f}\xb$,\index{F@$f$, the residue degree of $R$} the degree of $k$ over its prime subfield, is the \textit{residue degree} of $R$.\\

There must exist $m<n$ in $\xN$ with $\pi^m=\pi^n$. As $\pi^{n-m}$ is in the Jacobson radical of $R$, $\pi^{n-m}-1$ is a unit. So $\pi^m=0$. Let $\xh{a}$\index{A@$a$, the nilpotency index of the finite local PIR $R$} be the \textit{nilpotency index} of $R$, the smallest number with $\xp^{\xf[.25]a}=0$. For $0\ne r\in R$, we write $\boldsymbol{\ord_\xp}(r)=n$\index{ORD@$\ord_\xp$, ring element order with respect to $\xp$} in case $r\in\xp^{\xf[.25]n}-\xp^{\xf[.25]n+1}$. Then $r=\pi^ns$ for some $s\in R-\xp=\xU$. So if $0\subsetneq I\subseteq\xp$ is an ideal of $R$ and $0\ne r\in I$ an element with  minimal $\ord_p$, then $I=rR$. The only ideals of $R$ are therefore the $\pi^iR$ for $0\le i\le a$. We put $\ord_\xp(0)=a$. Let $\xh{e}\coloneqq\ord_\xp(p)$, the \textit{ramification index}\index{E@$e$, the ramification index of $R$} of $R$ over $\xZ$. Then the characteristic of $R$ is $p^{\xceil{a/e}}$, where $\ulcorner\xf\urcorner$ denotes the ceiling function. We select a unit $\xh{\varepsilon}$\index{YE@$\varepsilon$, a unit of $R$ for which $p=\varepsilon\pi^e$} for which $p=\varepsilon\pi^e$. (It is by no means unique, for elements such as $\varepsilon+\pi^{a-e}$ would work equally well.) Put $\xh{a_0}\coloneqq\xceil{a/e}$\index{A0@$a_0=\xceil{a/e}=$ the nilpotency index of the PIRs $R_0$ and $\xRo$}, and let $\xh{R_0}$\index{R0@$R_0$, the image of $\xZ$ in $R$} be the subring of $R$ generated by $1$. It is contained in all other subrings of $R$, and it is isomorphic to $\xZ/\xpe{a_0}\xZ$.\\

The residue map $R\twoheadrightarrow k$ is denoted by $\xh{\nu}$\index{YN@$\nu$, the residue map $R\to k$}. It is clear that $\nu\upharpoonright\xb\xU\xcol\xU\to\xu$ maps $\xU$ onto the group of units of $k$, which is cyclic of order $p^{\xf[.25]f}\xb-1$. The kernel $\xh{H_1}\coloneqq\ker(\nu\upharpoonright\xb\xU)=1+\,\xp$\index{H1@$H_1$, the group of one-units of $R$; its order is $p^{(a-1)f}$} consists of the so-called \textit{Einseinheiten}, aka \textit{one-units}. If $r=1+\pi r'\in H_1$, then $r^{\xf[.25]p}-1\in(\pi^2)$, since the binomial coefficient $\binom{p}{i}$ is a multiple of $p$ for $0<i<p$. Clearly, for some $n$, $r^{\xf[.25]p}\xexp{n}=1$, and the order $\xh{\mo}(r)$\index{O@$\mo$, order of a group element} of $r$ is a power of $p$. Hence so is $\xsize{H_1}$. Because $\xsize{\xu}$ is prime to $p$, the exact sequence $1\to H_1\to\xU\to\xu\to1$ splits. So $\xU$ is the product of $H_1$ and a cyclic subgroup $\xh{H_0}\cong\xu$\index{H0@$H_0$, a cyclic subgroup of $\xU$ isomorphic to $\xu$} of order $p^f\xb-1$. In particular, $\xsize{\xU}=(p^f\xb-1)p^{(a-1)f}$.\\

One has $H_0\subseteq\xh{\rho}\coloneqq\{r\in R\mid r^{\xf[.23]p^{\xb[1]\xf f}}\xb=r\}$\index{YR@$\rho$, the canonical system of representatives of $R$}. If $r\in\rho$ is a unit, its order divides $p^f\xb-1$, so $r$ must be in $H_0$. And if $r\notin\xU$ then $r\in\xp$, so $r$ is nilpotent. Since $r$ is equal to its $p^f$-th power, repeated application of that gives $r=0$. Hence $\rho=H_0\cup\{0\}$ and $\xsize{\rho}=\xsize{k}$. As $H_0\twoheadrightarrow\xu$, the map $\nu\upharpoonright\xb \rho\xcol\rho\to k$ is a bijection. $\rho$ is called the \textit{canonical system of representatives} of $R$ modulo $\xp$. Every $r\in R$ is congruent mod $\xp$ to a unique element of $\rho$. This feature holds in a more general setting, cf.\ Prop.\ II.8 in \cite{SE}.\\

For $n<a$, we have $\xp^n\xb[1.25]/\,\xp^{n+1}\xb=\pi^nR/\pi^{n+1}\xb R$, and the $\pi^nr$ for $r\in\rho$ represent the congruence classes of $\xp^n$ mod $\xp^{n+1}$. Thus the map $\mu_n\xcol k\to\xp^n\xb[1.25]/\,\xp^{n+1}$, with $\nu(r)\mapsto r\pi^n\xb$ mod $\xp^{n+1}\xb[.5]$, is well-defined and surjective. It is a homomorphism of $R$-modules. Since $k$ is a simple $R$-module, it follows that these maps $\mu_n$ are isomorphisms of $R$-modules.\\

Every element $r$ of $R$ can be uniquely written as a ``power series'' $r=\xsum{i=0}{a-1}r_i\pi^i$ where the $r_i$ run over $\rho$. Indeed, there is a unique $r_0\in\rho$ with $r-r_0\in\pi R$, say $r=r_0+\pi r'$, and we process $r'$ in the same way. As $\pi^{\xf[.5]a}=0$, this halts after $a$ such steps. To see uniqueness, if $\xsum{i=0}{a-1}r_i\pi^i=r=\xsum{i=0}{a-1}s_i\pi^i$ with the $s_i$ in $\rho$, then $r_0$ and $s_0$ both represent $r$, so they are equal and can be discarded. The congruence class of $\xsum{i=1}{a-1}r_i\pi^i=\xsum{i=1}{a-1}s_i\pi^i$ in $\xp/\xp^2\xb$ is represented by both $r_1\pi$ and $s_1\pi$, hence $r_1=s_1$. Repeating the argument, one finds $r_i=s_i$ for all $i$.\\

We will refer to this as \textit{uniqueness of representation}. As a direct consequence, $\xsize{R\xf}=\xpef{a}$. An $r=\xsum{i=0}{a-1}r_i\pi^i$ as above (with all $r_i\in\rho$) is a unit iff $r_0\ne0$.\\

For $n>1$, put $\xh{H_n}=\{\xf\delta\in \xU\mid\ord_\xp(\delta-1)\ge n\xf\}$\index{HN@$H_n$, the groups of higher one-units}, with $H_1$ part of the filtration $H_1\supseteq H_2\supseteq\cdots\supseteq H_a=\{1\}$. If $\delta=1+\pi^nr\in H_n$, with $n\ge1$ and $r\in R$, $\delta\mapsto\xo{r}$ defines a surjection $H_n\twoheadrightarrow k^+$ with kernel $H_{n+1}$. The $H_n$, with $n\ge2$, are known as the groups of the \textit{höhere Einseinheiten} (higher one-units). The similarity with local fields, whose terminology we borrow, is not accidental. We now summarize the basic properties of finite local PIRs.

\medskip
\begin{prop}\label{prop:basics}With the notation as above, the following hold.
\begin{enumerate}[label=\normalfont{(}\normalfont\arabic*)]
\item $R$ is of order $\xpef{a}$ and characteristic $\xpe{a_0}=p^{\xceil{a/e}}$. The only ideals of $R$ are the powers $\xp^i=\pi^iR$, with $0\le i\le a$.
\item The group of units $\xU$ is the direct product of a cyclic group $H_0$ of order $p^f\xb-1$ and a $p$-group $H_1$ of order $\xpef{(a-1)}$.
\item The canonical set of representatives $\rho$ of $R/\xp$ consists of the elements $r\in R$ that satisfy $r^{\xf[.23]p^{\xb[1]\xf f}}\xb=r$. One has $\rho=H_0\cup\{0\}$ and $\mid\xb[2.5]\rho\xb[2]\mid\,=p^f\xb = |k|$.
\item The $r\in R$ can be uniquely represented as $r=\xsum{i=0}{a-1}r_i\pi^i$ with the $r_i\in\rho$. Such an $r$ is a unit of $R$ iff $r_0\ne0$.
\item $0=\xp^a<\xp^{a-1}<\cdots<\xp^0=R$ is a composition series of the $R$-module $R$, with all composition factors $\xp^n/\xp^{n+1}\cong k$ as $R$-modules, for $0\le n<a$.
\item $1=H_a\lhd H_{a-1}\lhd\cdots\lhd H_1$ is a filtration of the group $H_1$, with elementary abelian factors $H_n/H_{n+1}\cong k^+$, for $1\le n<a$.
\item $\xsize{\xp^n}=\xsize{k}^{a-n}$ for $n<a$.
\item Every subring $R'$ of $R$ is a local ring, with maximal ideal $\xf\xp\cap R'$.
\item $R$ is $\xp$-adically complete.
\item $R$ is a chain ring: the set of its ideals is linearly ordered by inclusion.
\item $R$ is a field iff $a=1$. (In that case, $R$ coincides with $k$.)
\end{enumerate}
\begin{proof}
For (6), $p=0$ in $k$, so the additive group $k^+$ is clearly elementary abelian. (7) is immediate from (5). Alternatively, elements of $\xp^n$ can be uniquely written as $\xsum{i=n}{a-1}r_i\pi^i$, with the $r_i\in\rho$. To see (8), every subring $R'$ of $R$ is an Artin ring, hence a finite product of local Artin rings. But $R$ does not contain any non-trivial idempotents, so neither does $R'$. And $R'$ cannot have any prime ideals other than $\xp\xf[.7]\cap R'$. (9) and (10) are obvious; the ideals are even well-ordered. (11): $R$ is a field iff $\xp=0$ iff $\pi=0$ iff $a=1$. 
\end{proof}
\end{prop}

\medskip
The following is a direct application of Hensel's Lemma (\cite{HL}, in the formulation for general commutative rings).

\medskip
\begin{prop}\label{prop:hensel}
Let $t\in R[X]$ be a monic polynomial of degree $d$ such that $\xo{t}$ splits as a product of distinct monic irreducible polynomials in $k[X]$. Then $t$ has a unique factorization in $R[X]$ that lifts the one of $\xo{t}$.\\
If the factors are all linear, then $t$ has precisely $d$ zeroes in $R$.
\begin{proof}
Distinct monic irreducibles in $k[X]=(R/\xp)[X]$ are coprime, so the first claim follows immediately by Hensel's Lemma. Indeed, such factorizations lift uniquely to $(R/\xp^n)[X]$ for every $n\in\xN$, and taking $n=a$ yields a lifting to $(R/\xp^a)[X]=R[X]$. If all factors are linear, we have $t=\xprod{i=1}{d}(X-x_i)$ in $R[X]$. Then $\xo{t}=\xprod{i=1}{d}(X-\xo{x}_i)$ in $k[X]$, where $\xo{x}_i=\nu(x_i)$ for $1\le i\le d$. Clearly, the $x_i$ are distinct zeroes of $t$ in $R$. If $t(y)=0$ for some $y\in R$, we can factorize $t=(X-y)t_1$ for some $t_1\in R[X]$, by the Remainder Theorem. $\xo{y}$ must be one of the $\xo{x}_i$, say $\xo{y}=\xo{x}_1$. Since $t=(X-x_1)t_2$, where $t_2\coloneqq\xprod{j=2}{d}(X-x_j)$, by uniqueness of lifting $X-y=X-x_1$ and $t_1=t_2$. (In case $d=2$, the alternative $t_1=X-x_1$ and $X-y=t_2=X-x_2$ is excluded, for that would give $y=x_2$, so $\xo{x}_1=\xo{y}=\xo{x}_2$). Thus $y=x_1$.
\end{proof}
\end{prop}

\medskip
The lemma below generalizes a well-known fact from elementary number theory.

\medskip
\begin{lem}\label{lem:elem}
Let $r\in R-\{1\}$ have $\ord_\xp(r-1)=m$. If $mp>m+e$, then for every $n\in\xN$ one has $r^{p^n}\xb=1$ or $\ord_\xp(r^{p^n}\xb-1)=m+ne$, or both.
\begin{proof}
$r=1+\delta\pi^m$ for some $\delta\in \xU$. Because $p=\varepsilon\pi^e$, $r^p=(1+\delta\pi^m)^p\equiv 1+\delta\xf[.25]\varepsilon\pi^{m+e}$ mod $\pi^{m+e+1}$\xb, for $p$ divides the non-trivial binomial coefficients, and $mp\ge m+e+1$. But $\delta\xf[.25]\varepsilon\in\xU$. This settles the case $n=1$. When it holds for $n-1$, by $(m+(n-1)e)\xf[.25]p>m+(n-1)e+e$ the property follows for $n$ too.
\end{proof}
\end{lem}

\begin{xmp}\label{xmp:1}Let $m=\xceil{(e+1)/(p-1)}$ and $r=1-\pi^m\in H_1$. Then $mp>m+e$, so by Lemma \ref{lem:elem} either $r^{p^n}\xb=1$ or $\ord_\xp(r^{p^n}\xb-1)=m+ne<a$, for every $n\in\xN$. Hence the order of $r$ in $\xU$ is $p^n$, where $n=\xceil{(a-m)/e}$. If $p\ge e+2$, one has $m=1$. So in that case $r\coloneqq1-\pi$ and its inverse $\xsum{i=0}{a-1}\pi^i$ are of the indicated order in $\xU$.
\end{xmp}

\section{Ring structure}\label{sec:structure}

The purpose of this section is to construct algebraic number fields $K\subseteq L$ and surjections from their rings of integers $\xO_K$ and $\xO_L$ to certain subrings of $R$:\\

\begin{adjustbox}
{center}\begin{tikzpicture}[scale=1,text height=1.5ex, text depth=0.25ex, >= Straight Barb]
\node (OQ) at (1,2.2) {$\xZ$}; 
\node (OK) at (2.5,2.2) {$\xO_K$};
\node (OL) at (4,2.2) {$\xO_L$};
\node (R0) at (1,1) {$R_0$};
\node (Ro) at (2.5,1) {$\xRo$};
\node (R) at (4,1) {$R$};
\path[-stealth, auto] (OQ) edge[subseteq] (OK);
\path[-stealth, auto] (OK) edge[subseteq] (OL);
\path[-stealth, auto] (R0) edge[subseteq] (Ro);
\path[-stealth, auto] (Ro) edge[subseteq] (R);
\draw[->>]
(OQ) edge (R0)
(OK) edge node[right,pos=0.4] {$\varphi$} (Ro)
(OL) edge node[right,pos=0.4] {$\psi$} (R);
\end{tikzpicture}\end{adjustbox}

\medskip
The degree $[K:\xQ]$ will be $f$, with $p$ inert in $K$. So $p\xf[.5]\xO_K$ will have residue degree $f$ over $p\xZ$. And $[L:K]$ will be equal to $e$, with $p\xf[.5]\xO_K$ totally ramified in $L$. The intermediate ring $\xh{\xRo}\coloneqq R_0[\xf[.5]\rho\xf[.5]]$\index{Ro@$R_\omega$, the coefficient ring of $R$} is the subring of $R$ generated by the canonical representatives.\\

Let $\xh{h}\in\xZ[X]$\index{HZ@$h$, a degree $f$ monic in $\xZ[X]$ that is irreducible in $\xF_p[X]$} be any monic of degree $f$ that is irreducible in $\xF_p[X]$. Put $\xh{K}\coloneqq\xQ[X]/(h)$\index{K@$K$, a field extension of degree $f$ of $\xQ$ in which $p$ is inert}, with ring of integers $\xh{\xO_K}$\index{OK@$\xO_K$, the ring of integers of $K$}, and let $\xh{\beta}\in\xO_K$\index{YB@$\beta$, a zero of the polynomial $h$ in $\xO_K$} be a zero of $h$. By the Kummer-Dedekind theorem (Th.\ 8.2 in \cite{ST}, or Th.\ \ref{appt:KD} in the appendix here), $p\xf[.3]\xZ[\beta]$ is a prime ideal of the order $\xZ[\beta]$ of $\xO_K$, since $h$ is irreducible mod $p$. Prop.\ \ref{appp:prinv} then shows that $\xh{\xP_K}\coloneqq p\xf[.3]\xO_K$\index{PK@$\xP_K$, the unique $\xO_K$-prime that lies over $p$} is a prime ideal of $\xO_K$, i.e., $p$ is \textit{inert} in $K$, that $\xZ[\beta]_{p\xZ[\beta]}$ $=$ $\xO_{K,{\xP_K}}$ is a DVR, and that $\xZ[\beta]/p\xZ[\beta]$ $\cong$ $\xO_K/\xP_K$. The latter residue field is isomorphic to $\xF_{p^f}$\xb, i.e.\ to $k$. Identifying these fields, the element of $k$ that corresponds to $\beta$ mod $\xP_K$ will be denoted by $\xs{\xh{\beta}}$.\index{YB1@$\xs{\beta}$, the zero in $k$ of the polynomial $h$ that corresponds to $\beta\in\xO_K$}\\

In $k[X]$, $h$ splits into $f$ distinct monic linear factors. By Prop.\ \ref{prop:hensel}, the same is true in $R[X]$, and $h$ has precisely $f$ zeroes in $R$. Let $\xh{b}\in R$\index{B@$b$, the zero in $R$ of the polynomial $h$ that corresponds to $\beta\in\xO_K$} be the zero that maps to $\xo{\beta}\in k$. Since $\xZ[\beta]\cong\xZ[X]/(h)$ and $h(b)=0$, putting $\beta\mapsto b$ gives a well-defined ring homomorphism $\xh{\varphi}\xcol\xZ[\beta]\to R$. The image, $\xRb$, is a local ring by Prop.\ \ref{prop:basics}.8. As $\xRb/(p)$ is a quotient of $\xZ[X]/(p,h)\cong\xF_{p^f}\xb$, a field, $p\xRb$ is a prime ideal of $\xRb$. By the Lemma below, $\xRb$ is a finite local PIR. As $p$ generates its maximal ideal, $p$ is a uniformizer for $\xRb$.

\medskip
\begin{lem}\label{lem:principir}
An Artinian local ring $A$ with principal maximal ideal $\mathfrak{m}$ is a PIR.
\begin{proof}
Let $\mathfrak{m}=xA$. As $A$ is Artinian, there exists an $n$ such that $\mathfrak{m}^{n+1}=\mathfrak{m}^n$. By Nakayama's lemma, $\mathfrak{m}^n=0$. Let $I$ be a non-zero ideal of $A$. Take $i\ge0$ maximal with respect to $I\subseteq\mathfrak{m}^i$, and take $y\in I$ with $y\notin\mathfrak{m}^{i+1}$. Then $y=x^iz$ for some $z\in A-\mathfrak{m}=A^\times\xb$, and therefore $\mathfrak{m}^i=x^iA=yA\subseteq I$.\end{proof}
\end{lem}

The attributes $(a,e,f,p)$ that pertain to $R$, forming its \textit{type}, take the values $(a_0,1,f,p)$ for the finite local PIR $\xRb$. In particular, one has $|\xRb|=\xpef{a_0}$. Prop.\ \ref{prop:basics}.3, applied to $\xRb$, tells us that $\rho\subseteq \xRb$, hence $\xRo\subseteq \xRb$. Since the $\xpef{a_0}$ sums $\xsum{i=0}{a_0-1}r_ip^i$ (with the $r_i\in\rho$) are all different in $R$ and belong to $\xRo$, it follows that $\xRo=\xRb$. So $\xRb$ does not not depend on the choice of the particular zero $b$ of $h$ in $R$, and all $f$ zeroes of $h$ in $R$ are already in $\xRo$. As it has ramification index 1, the finite local PIR $\xRo$ is called \textit{unramified}.\\

$\xRo$ is a \textit{coefficient ring} of the complete local ring $R$ (\cite{CR}): it is a complete local ring contained in $R$, with the same residue field as $R$, and with maximal ideal $\xRo\cap\xf\xp=p\xRo$, where $p$ is the characteristic of the residue field. We shall refer to $\xRo$ as \textit{the} coefficient ring. Indeed, every coefficient ring $C\subseteq R$ is a PIR by Lemma \ref{lem:principir}, for its prime ideal $p\xf[.5]C$ is principal. Because $C/p\xf[.5]C=k$, $\rho\subseteq C$ (by Prop.\ \ref{prop:basics}.3), hence $\xRo\subseteq C$. But $C$ must be of type $(a_0,1,f,p)$, and hence of size $\xpef{a_0}=|\xRo|$. It follows that $C=\xRo$.\\

We have $\xO_K\subseteq\xO_{K,\xP_K}=\xZ[\beta]_{p\xZ[\beta]}\overset{\varphi}\longrightarrow\xRo$, producing the homomorphism $\xh{\varphi}\xcol\xO_K\to\xRo$\index{YX@$\varphi$, the structure homomorphism $\xO_K\twoheadrightarrow\xRo$} in the diagram above. It is surjective because $\varphi\xcol\xZ[\beta]\to\xRo$ is.\\

\textbf{Note}. If $V=\widehat{\xO_K}_{\xP_K}$ is the completion of $\xO_K$ at $\xP_K$, a $v$-ring, the map $\varphi$ can be extended to $V\twoheadrightarrow\xRo$. This effectively copies the proof of Cohen's structure theorem in the unequal-characteristic case, Th.\ 11 of \cite{CO}, applied to the complete local ring $R$.

\medskip
\begin{prop}\label{prop:omega}With notation as above, $R$ has a unique coefficient ring, $\xRo=R_0[\xf[.5]\rho\xf[.5]]=\xRb$, having the following properties.
\begin{enumerate}[label=\normalfont{(}\normalfont\arabic*)]
\item $\xRo$ is an unramified finite local PIR of type $(a_0,1,f,p)$ with uniformizer $p$ and residue field $k$.
\item $R$ and $\xRo$ have the same canonical set of representatives $\rho$.
\item $R=\xRo[\pi]$.
\item $R=\xRo$ iff $e=1$.
\item $\xRo\cong\xZ[\beta]/(\xpe{a_0})\cong\xO_K/(\xpe{a_0})$, with $\xZ[\beta]\twoheadrightarrow\xRo$ given by $\beta\mapsto b$.
\item $\xRo$ is determined up to isomorphism by $a_0$, $f$ and $p$.
\end{enumerate}
\begin{proof}
Item (3) is immediate from (2) and point (4) of Prop.\ \ref{prop:basics}. For (4): if $R=\xRo$ then $p$ is a uniformizer for $R$. As $p\sim\pi^e$, it follows that $e=1$. Conversely, if $e=1$ then $a_0=a$, so $\xsize{\xRo}=\xsize{R}$. Since the additive group of $\xZ[\beta]$ is free of rank $f$, we have $[\xZ[\beta]:\xpe{a_0}\xZ[\beta]]=\xpef{a_0}=\xsize{\xRo}=\mathrm{N}_{K/\xQ}(\xpe{a_0})=[\xO_K:\xpe{a_0}\xO_K]$. Since $\xpe{a_0}\in\ker(\varphi\xb[1]\upharpoonright\xb\xZ[\beta])$, (5) results. And (6) follows directly from (5).
\end{proof}
\end{prop}

\medskip
Let $\xh{\vartheta}\xcol\xRo[Y]\twoheadrightarrow R$ be the homomorphism defined by $Y\mapsto\pi$. When $R$ is \textit{equal-characteristic}, that is, $\xchar(R)=p=\xchar(k)$, that is, $a=e$, by Prop.\ \ref{prop:omega}.1 and Prop.\ \ref{prop:basics}.11 we have $\xRo=k$. Clearly, $\ker(\vartheta)=(Y^a\xb)$, hence $R\cong k[Y]/(Y^a\xb)$.

\medskip
\begin{cor}\label{cor:eqchar}In the equal-characteristic case $a=e$ one has:
\begin{enumerate}[label=\normalfont{(}\normalfont\arabic*)]
\item$R$ is isomorphic to $k[Y]_{(Y)}/(Y^a\xb)$, a quotient of the discrete valuation ring $k[Y]_{(Y)}$.
\item The structure of $R$ is determined by $a$, $f$ and $p$.$\hfill\square$
\end{enumerate}
\end{cor}

\medskip
We now concentrate on the case $\xh{a>e}$.\\

$\xO_K\overset{\varphi}\longrightarrow\xRo$ and $\xRo[Y]\overset{\vartheta}\longrightarrow R$ give rise to a homomorphism $\xh{\chi}\xcol\xO_K[Y]\to R$\index{YY@$\chi$, a homomorphism $\xO_K[Y]\to R$ with $\beta\mapsto b$ and $Y\mapsto\pi$}, which is defined on $\xZ[\beta][Y]$ by $\beta\mapsto b$ and $Y\mapsto\pi$. The image of $\chi$ is a subring of $R$ that contains both $\xRb=\xRo$ and $\pi$. Thus $\chi$ is surjective.\\

Using Prop.\ \ref{prop:basics}.4, write $\varepsilon=\xsum{i=0}{a-1}r_i\pi^i$ with the $r_i\in\rho$ and $r_0\ne0$. Put $\xh{\tilde{\varepsilon}(Y)}\coloneqq\xsum{i=0}{a-e-1}r_i Y^i\in\xRo[Y]$. Then $\deg(\tilde{\varepsilon})<a-e$, $\tilde{\varepsilon}(0)=r_0$ is in $\xUo$, and $\vartheta(\tilde{\varepsilon}(Y)Y^e)=\xsum{i=0}{a-e-1}r_i\pi^{i+e}=\varepsilon\pi^e=p$. So $p-\tilde{\varepsilon}(Y)Y^e\in\ker(\vartheta)$. Clearly, $Y^a\xb\in\ker(\vartheta)$ as well.

\medskip
\begin{prop}\label{prop:kernel}
Let $a>e$. With $\tilde{\varepsilon}$ as above, $\ker(\vartheta)=(Y^a,p-\tilde{\varepsilon}(Y)Y^e)$.
\begin{proof}
Let $t\in\ker(\vartheta)$. As $\pi^a=0$, we can assume $\deg(t)<a$, say $t=\xsum{j=0}{a-1}c_jY^j$, with the coefficients $c_j\in\xRo$. Because $t(\pi)=0$, we have $c_0\in\xp\cap\xRo=p\xRo$, so $c_0=\delta p^i$ for some $\delta\in\xUo$ and $i\in\xN$. We then replace the monomial $c_0$ in $t$ with $\delta\xf[.75]\tilde{\varepsilon}(Y)^iY^{ei}$, reduce the result mod $Y^a\xb$, and denote the coefficients by $c_j$ again. Then $t$ is still in $\ker(\vartheta)$, but now $c_0=0$. For $j\ge1$, we have $c_j=r+s$ for a unique $r\in\rho$ and an $s\in p\xRo$. Write $s=\delta p^i$ in $\xRo$, with $\delta\in\xUo$ and $i\in\xN$. Replace the monomial $c_jY^j$ by $rY^j+\delta\xf[.75]\tilde{\varepsilon}(Y)^iY^{ei+j}$, and reduce mod $Y^a$ again. This does not introduce new terms of degree lower than $j$, and the coefficient of $Y^j$ is now in $\rho$. After $a$ steps, all coefficients of $t$ are in $\rho$, and $t$ is still of degree less than $a$. By unique representation, $t$ must be zero, and the claim follows.
\end{proof}
\end{prop}

\medskip
Consequently, $R\cong\xRo[Y]/(Y^a\xb,p-\tilde{\varepsilon}(Y)Y^e)=\xRo[Y]/(Y^a,g_0)$, where $\xh{g_0(Y)}\coloneqq Y^a+\tilde{\varepsilon}(Y)Y^e-p$ is monic of degree $a$ over $\xRo$. Modulo $p$, $g_0\equiv Y^e\cdot(Y^{a-e}+\tilde{\varepsilon}(Y))$. Since $\tilde{\varepsilon}(0)$ is a unit, the two factors are coprime. As $(p)\in\max(\xRo)$, by Hensel's Lemma this factorization lifts (uniquely) to $\xRo[Y]$ as $g_0=g\cdot g_1$, where both factors are monic, and the first one, $\xh{g(Y)}$, is of degree $e$ and satisfies $g(Y)\equiv Y^e$\index{G@$g$, a $p\,$-Eisenstein polynomial over $\xRo$ of degree $e$} mod $p$ and $g(0)\sim p$ in $\xRo$.\\

As a polynomial in $\xRo[Y]$, $g$ is $p\xf[2]$-\textit{Eisenstein} (of degree $e$) in the usual sense: $g\equiv Y^e$ mod $p$ and $g(0)\in (p)-(p^2)$. Thus $g\in\xh{\nabla}\coloneqq\{q(Y)\in\xRo[Y]\mid \deg(q)=e$ and $q$ is $p\,$-Eisenstein$\}$\index{ZN@$\nabla$, the set of $p\,$-Eisenstein polynomials of degree $e$ over $\xRo$}.\\

In $R$, we now have $0=g_0(\pi)=g(\pi)g_1(\pi)$. The second factor is a unit, because $g_1\equiv Y^{a-e}+\tilde{\varepsilon}(Y)$ mod $p$. It follows that $g(\pi)=0$ in $R$.\\

Let $\xh{R_g}\coloneqq\xRo[Y]/(g)$\index{RG@$R_g$, the neat PIR associated to the Eisenstein polynomial $g$}. Then the homomorphism $R_g\to R$, $\xs{Y}\mapsto\pi$, is well-defined and surjective. If $\mathfrak{m}$ is a prime ideal of $R_g$, then $p\in\mathfrak{m}$, so $\xs{Y}\in\mathfrak{m}$, for $g\equiv Y^e$ mod $p$. But since $g(0)\sim p$, we have $R_g/(\xs{Y})\cong\xRo/(p)\cong k$, using Prop.\ \ref{prop:omega}.1. Therefore, $\mathfrak{m}=(\xs{Y})$. Thus $R_g$ is a local PIR with uniformizer $\xs{Y}$ and residue field $k$. As $g$ is monic, $\xsize{R_g}=\xsize{\xRo}{}^e=\xpe{a_0\xb[.3]fe}$. Due to Prop.\ \ref{prop:basics}.1, the nilpotency index of $R_g$ is $a_0e$. Writing $g=Y^e+\xsum{i=0}{e-1}p\xf[.3]c_iY^i$, with the $c_i\in\xRo$ and $c_0\in\xUo$, one has $\xs{Y}{}^e=-p\cdot(\xsum{i=0}{e-1}c_i\xs{Y}{}^i)\sim p$ in $R_g$. Hence $R_g$ is of ramification index $e$. The kernel of $R_g\twoheadrightarrow R$ is $(\xs{Y}{}^n)$ for some $n$. As $\pi^{a-1}\xb\ne0$ in $R$, necessarily $n=a$.

\medskip
\begin{prop}\label{prop:Rg}Let $a>e$. Then $g(\pi)=0$ in $R$. Moreover:
\begin{enumerate}[label=\normalfont{(}\normalfont\arabic*)]
\item$R_g=\xRo[Y]/(g)$ is a finite local PIR of type $(a_0e,e,f,p)$, of which $R$ is a quotient: $R\cong R_g/(\xs{Y}{}^a)$.
\item Conversely, for any $q\in\nabla$, the ring $\xh{R_q}\coloneqq\xRo[Y]/(q)$ is a finite local PIR of type $(a_0e,e,f,p)$ that has $\xs{Y}$ as a uniformizer. Its quotient $R'\coloneqq R_q/(\xs{Y}{}^a)$ is a local PIR of type $(a,e,f,p)$.
\end{enumerate}
\begin{proof}
(2) follows as above, when one replaces $g$ by $q$.
\end{proof}
\end{prop}

\medskip
Let $\xh{\tilde{g}}\in\xZ[\beta][Y]$ be a monic lift of degree $e$ of $g\in\xRo[Y]\cong\xZ[\beta]/(\xpe{a_0})[Y]$. Since $g$ is $p\,$-Eisenstein and $a_0\ge2$, it follows that $\tilde{g}\in\xO_K[Y]$ is $p\,$-Eisenstein too. For $p\xf[.5]\xO_K$ is a maximal ideal and $\tilde{g}(0)\in p\xf[.5]\xO_K-p^2\xO_K$. Hence $\tilde{g}$ is irreducible in $\xO_K[Y]$, so, since $\xO_K$ is integrally closed, also in $K[Y]$, by Gauss' Lemma (cf.\ Prop.\ \ref{appp:GL}). The map $\chi\xcol\xO_K[Y]\to R$ factors through $\xO_K[Y]\to\xO_K[Y]/(\tilde{g})$ as $g(\pi)=0$. Put $\gamma\coloneqq Y$ mod $(\tilde{g})$ and let $\xh{L}\coloneqq K[Y]/(\tilde{g})=K(\gamma)$\index{L@$L$, an extension of $K$ of degree $e$ in which $p\xf[.5]\xO_K$ is totally ramified} be the quotient field of $\xO_K[Y]/(\tilde{g})$, with ring of integers $\xh{\xO_L}$\index{OL@$\xO_L$, the ring of integers of $L$}. We obtain (the solid part of):\\

\begin{adjustbox}
{center}\begin{tikzpicture}[scale=1,text height=1.5ex, text depth=0.25ex, >= Straight Barb]
\node (OKY) at (1,3) {$\xO_K[Y]$}; 
\node (R) at (2.8,3) {$R$};
\node (Q) at (1,1.73) {$\xO_K[Y]/(\tilde{g})$};
\node (OKG) at (2.8,1.73) {$\xO_K[\gamma]$};
\node (OL) at (4.5,1.73) {$\xO_L$};
\node (OO) at (4.5,1.93) {};
\draw[->] (OKY) edge (Q);
\draw[->>]
(OKY) edge node[above,pos=0.4] {$\chi$} (R)
(OKG) edge node[left,pos=0.4] {$\psi\xb[2]$} (R);
\draw[->] (OO) edge[dashed] node[left,pos=0.34] {$\psi\xf[5]$} (R);
\path (Q)--(OKG) node[midway]{=};
\path[-stealth, auto] (OKG) edge[subseteq] (OL);
\end{tikzpicture}
\end{adjustbox}

\vspace{3pt}
\medskip
Then $\psi\xcol\xO_K[\gamma]\twoheadrightarrow R$ extends $\varphi\xcol\xO_K\to\xRo$, sends $\beta\mapsto b$ and $\gamma\mapsto\pi$, and it extends to $\xO_L$ as follows.\\

We claim that the order $\xO_K[\gamma]$ of $\xO_L$ has a unique prime $\xP_\gamma\coloneqq p\xf[.5]\xO_K[\gamma]+\gamma\xf[.5]\xO_K[\gamma]$ lying over $\xP_K$, this prime is an invertible ideal of the order, and one has $p\xf[.5]\xO_K[\gamma]=\xP_\gamma^{\xf[.3]e}$. This follows from the Kummer-Dedekind theorem (\ref{appt:KD}). Indeed, the image of $\tilde{g}$ in $(\xO_K/\xP_K)[Y]\cong k[Y]$ is $Y^e$. And since $\tilde{g}$ is $p\,$-Eisenstein, the remainder on dividing $\tilde{g}$ by $Y$ in $\xO_K[Y]$ is not in $p^2\xO_K=\xP_K^{\xf[.2]2}$. Now let $\xh{\xP_L}\coloneqq\xP_\gamma\xO_L=p\xf[.5]\xO_L+\gamma\xf[.5]\xO_L$\index{PL@$\xP_L$, the unique $\xO_L$-prime that lies over $p$}. Then by Prop.\ \ref{appp:prinv}, $p\xf[.5]\xO_L=\xP_L^{\xf[.3]e}$ is the prime factorization of $p$ in $\xO_L$, and $\xO_K[\gamma]/\xP_\gamma\cong\xO_L/\xP_L$, and we have $\xO_K[\gamma]_{\xP_\gamma}=\xO_{L,\xP_L}\xb$, a DVR. Localizing extends $\psi\xcol\xO_K[\gamma]\to R$ to $\psi\xcol\xO_K[\gamma]_{\xP_\gamma}\to R$, which in turn restricts to the desired $\xh{\psi}\xcol\xO_L\twoheadrightarrow R$\index{YZ@$\psi$, the structure homomorphism $\xO_L\twoheadrightarrow R$}, in view of $\xO_L\subseteq\xO_{L,\xP_L}=\xO_K[\gamma]_{\xP_\gamma}\xb$.\\

Thus $\xP_K$ is totally ramified in $\xO_L$ with index $e$. The residue degree and ramification index of $\xP_L$ over $p\xZ$ are $f$ and $e$, respectively. And $[L:\xQ]=ef$. Note that, for $0\le i\le a$, $\psi(\gamma^i\xpe{a-i})=\pi^i\xpe{a-i}\sim\pi^{i+(a-i)e}=0$, since $i+(a-i)e\ge a$. So $\xP_\gamma^{\xf[.3]a}\subseteq\ker(\psi)$. But $\xsize{\xO_K[\gamma]/\xP_\gamma^{\xf[.3]a}}=[\xO_K[\gamma]:\xP_\gamma]{}^a$ (Prop.\ \ref{appp:absnorm}) $=[\xO_L:\xP_L]{}^a$ $=\mathrm{N}_{L/\xQ}(\xP_L)^a=\xpef{a}=\xsize{R}$. Thus $\ker(\psi\xf[.7]{\upharpoonright}\xf[.7]\xO_K[\gamma])=\xP_\gamma^{\xf[.3]a}$. As $\tilde{g}$ is Eisenstein, one has $p\in\gamma\xO_{L,\xP_L}$, so $\xP_L\xO_{L,\xP_L}=\gamma\xO_{L,\xP_L}$ and $\ker(\xO_{L,\xP_L}\underset{\psi}\longrightarrow R)=\gamma^a\xO_{L,\xP_L}\xb$.\\

We summarize.

\medskip
\begin{thm}\label{thm:main}For a finite local PIR $R$ one has, with notation as above:
\begin{enumerate}[label=\normalfont{(}\normalfont\arabic*)]
\item $R$ is the homomorphic image of a discrete valuation ring $D$ of an algebraic number field. In effect, there exist field extensions $\xQ\subseteq K=\xQ[X]/(h(X))\subseteq L=K[Y]/(\tilde{g}(Y))$, such that $[K:\xQ]=f$, $p$ is inert in $\xO_K$, $[L:K]=e$, and $p\xf[.5]\xO_K$ is totally ramified in $\xO_L$. For $D$ one can take the DVR $\xO_{L,\xP_L}\xb$, where $\xP_L$ is the unique prime of $\xO_L$ lying over $p$. One has $R\cong D/\xP_L^{\xf[.5]a}D=D/(\xs{Y}{}^a)$. The ramification index and residue degree of the extension $L/\xQ$ with respect to $p$ match those of the PIR $R$.
\item For every monic polynomial $\hat{h}\in\xZ[X]$ of degree $f$ and irreducible mod $p$ there exists an isomorphism $\xZ[X,Y]/(\xpe{a_0},Y^a\xb,\hat{g}(X,Y),\hat{h}(X))\xrightarrow{\sim}R$ for a suitable polynomial $\hat{g}(X,Y)$, which is monic and of degree $e$ with respect to $Y$ and has a $p\,$-Eisenstein image in $(\xZ[X]/(\hat{h}))[Y]$. Given a zero $\hat{b}$ of $\hat{h}$ in $R$, the isomorphism is defined by $\xs{X}\mapsto\hat{b}$ and $\xs{Y}\mapsto\pi$.
\end{enumerate}
\begin{proof}
(1) has just been established in the case $a>e$. The construction copies easily to the case $a=e$. Indeed, just take $\tilde{g}\coloneqq Y^e-p\in\xZ[\beta][Y]$. It is irreducible in $K[Y]$. We again put $L\coloneqq K[Y]/(\tilde{g})$, construct the structure map $\psi\xcol\xO_L\twoheadrightarrow R$ that sends $\beta\mapsto b$ and $\gamma(\coloneqq Y$ mod $(\tilde{g}))\mapsto\pi$, and argue as before. We define $g$ to be the image of $\tilde{g}$ in $\xRo[Y]$ here. As $p=0$ in $R$, we have $g=Y^e\xb=Y^a\xb$, and hence $R\cong k[Y]/(Y^a\xb)$ (Cor.\ \ref{cor:eqchar}) $=\xRo[Y]/(Y^a\xb)=\xRo[Y]/(Y^a\xb,g)$.
When $a>e$, $R\cong\xRo[Y]/(Y^a\xb,g)$ holds as well, by Prop.\ \ref{prop:Rg}. So in the general case ($a\ge e$) one has $R\cong\xRo[Y](Y^a\xb,g)\cong\xZ[\beta,Y]/(\xpe{a_0},Y^a\xb,\tilde{g})$ (cf.\ Prop.\ \ref{prop:omega}.5). Now notice that the only properties of the polynomial $h\in\xZ[X]$ we have used are that it is monic of degree $f$ and irreducible mod $p$. So for (2) we may assume that $\hat{h}=h$. Then $\xZ[\beta]=\xZ[X]/(\hat{h}(X))$. If $\hat{g}(X,Y)\in\xZ[X,Y]$ is a lift of $g(Y)\in\xRo[Y]$ under $\xZ[X,Y]\twoheadrightarrow\xRo[Y]$, $X\mapsto\beta$, statement (2) follows.
\end{proof}
\end{thm}

\medskip
Below, we again restrict to the unequal-characteristic case $\xh{a>e}$. (Note that when $a=e$, one has $p=0$ in $\xRo$ and the notion of $p\,$-Eisenstein polynomials in $\xRo[Y]$ does not make much sense.)\\

By Prop.\ \ref{prop:omega}.6, $\nabla$ depends only on the attributes $(a,e,f,p)$ of $R$. Therefore, $\xh{\Delta}\coloneqq\{\,q(Y)\in\nabla\mid q$ has a zero in $R\,\}$\index{YD@$\Delta$, the polynomials in $\nabla$ that have a zero in $R$} is a well-defined invariant of $R$. Since $g(\pi)=0$ by Prop.\ \ref{prop:Rg}, we have $g\in\Delta$.

\medskip
\begin{prop}\label{prop:delta}Let $a>e$. If $q\in\Delta$ and $u\in R$ is a zero of $q$, then $u$ is a uniformizer of $R$.
\begin{proof}
Let $u\in R$ and $q=\xsum{i=0}{e-1}p\xf[.3]d_iY^i+Y^e\in\Delta$, where all $d_i\in\xRo$ and $d_0\in\xUo$, such that $q(u)=0$. Then $u^e\in pR$, and hence $\ord_\xp(u)>0$. Thus $e\cdot\ord_\xp(u)=\ord_\xp(p\xf[.3]d_0+\xsum{i=1}{e-1}p\xf[.3]d_iu^i)=\ord_\xp(p\xf[.3]d_0)=e$, and therefore $\ord_\xp(u)=1$. That is, $u$ is a uniformizer of $R$.
\end{proof}
\end{prop}

\medskip
This $\Delta$, along with the attributes $(a,e,f,p)$ we defined earlier, makes up $R$'s \textit{signature} $(a,e,f,p,\Delta)$. We will study $\Delta$ in more detail in section \ref{sec:delta}.

\medskip
\begin{thm}\label{thm:main2}If $a>e$, the finite local PIR $R$ is determined up to isomorphism by its signature $(a,e,f,p,\Delta)$.
\begin{proof}
$\xRo$ is known because of Prop.\ \ref{prop:omega}.6. Take any $q\in\Delta$ and $u\in R$ such that $q(u)=0$. By Prop.\ \ref{prop:delta}, $u$ is a uniformizer of $R$. So $R=\xRo[u]$, by Prop.\ \ref{prop:omega}.3. Then $R_q/(\xs{Y}{}^a)\to R$, $\xs{Y}\mapsto u$, is well-defined and surjective. For cardinality reasons (cf.\ Prop.\ \ref{prop:basics}.1 and Prop.\ \ref{prop:Rg}.2), it is injective as well.
\end{proof}
\end{thm}

\medskip
We discuss a few examples.

\begin{xmp}\label{xmp:2}$R\coloneqq\xZ[X]/(8,X^2-2)$ is a type $(6, 2, 1, 2)$ finite local PIR. This is one of the two rings Mr.\ Page supplied as examples to point out that $a$, $e$, $f$ and $p$ alone do not determine $R$. We have $a_0=3$, and we can take $\pi=\xs{X}$ and $\varepsilon=1$ mod 8. Because $f=1$, $\xRo=\xZ/8\xZ$ and we can put $K=\xQ$, $h=X-1$, $\beta=1$ and $b=1$ mod 8. We can use the constant $1$ mod 8 for the polynomial $\tilde{\varepsilon}(Y)$ here. Then $g_0=Y^6+Y^2-2$, which factorizes in $(\xZ/8\xZ)[Y]$ as $(Y^2-2)(Y^4+2Y^2-3)$. So $g=Y^2-2$, $\gamma=\sqrt{2}$, $L=\xQ(\gamma)$ and $\xO_L=\xZ[\sqrt{2}]$. The homomorphism $\psi\xcol\xO_L\to R$, with $\beta\mapsto b$ and $\gamma\mapsto\pi$, is a surjection with kernel generated by $\gamma^a=8$. We have $R\cong R_g/(\xs{Y}{}^a)=(\xZ/8\xZ)[Y]/(Y^2-2,Y^6\xb)=(\xZ/8\xZ)[Y]/(Y^2-2)$. The latter equality is explained by Prop.\ \ref{prop:neat}.1 below.
\end{xmp}

\begin{xmp}\label{xmp:3}$R\coloneqq\xZ[X]/(8,X^2+2)$, Mr.\ Page's other example, is a $(6, 2, 1, 2)$ as well. Much the same parameters apply, except that $\varepsilon=-1$ mod 8, $g(Y)=Y^2+2$ and $g_1=Y^4-2Y^2+3$. Here, $\psi$ is the map $\xZ[\sqrt{-2}]\to R$ with $\gamma=\sqrt{-2}\mapsto\xs{X}=\pi$ and kernel $(\gamma^6)=(-8)$. For this ring, $q(Y)\coloneqq Y^2-2$ is an Eisenstein polynomial of degree $e$ in $\xRo[Y]$ that has no zeroes in $R$. I.e., $q\in\nabla-\Delta$.
\end{xmp}

\begin{xmp}\label{xmp:4}Slightly less trivially, if $p=f=3$ we can take $h=X^3-X-1$, as it is irreducible mod 3. Assuming $a=3$, $e=2$ and $\varepsilon=1$ in $R$, we have $a_0=2$ and $g_0=Y^3+Y^2-3=(Y^2+3Y-3)(Y-2)$ in $(\xZ/9\xZ)[Y]$. This tells us that $\xRo=\xZ[X]/(9,X^3-X-1)$ and $R=\xRo[Y]/(Y^3,Y^2+3Y-3)$.
\end{xmp}

We will revisit these examples in the next section.

\section{The invariant \texorpdfstring{$\Delta$}{Delta}}\label{sec:delta}

There exists a natural group action on the set $\nabla$ of $p\,$-Eisenstein polynomials over $\xRo$ of degree $e$ under which $\Delta$ is a union of orbits. The finite group involved is non-commutative in general. It depends only on $(a,e,f,p)$. In case $e\mid a$, $\Delta$ coincides with a single orbit, the one of $g\in\Delta$. But this may also hold when that condition does not apply. We again assume $\xh{a>e}$. Some preparation is in order before we can present the group action. $\xh{U}\coloneqq\{u\in R\mid\ord_\xp(u)=1\}$\index{U@$U$, the set of uniformizers of $R$} denotes the set of uniformizers of $R$.\\

Let $q\in\Delta$. I.e., $q\in\xRo[Y]$ is a $p\,$-Eisenstein polynomial of degree $e$ with a zero $u\in R$. By Prop.\ \ref{prop:delta}, we have $u\in U$. Conversely, just like $\pi$ is a zero of the Eisenstein $g\in\Delta$, every uniformizer $u$ is a zero of some $q\in\Delta$. We then have a homomorphism $R_q=\xRo[Y]/(q)\twoheadrightarrow R$ given by $\xs{Y}\mapsto u$. The kernel is $(\xs{Y}{}^a)$ and $R_q$ is an $(a_0e,e,f,p)$ by Prop.\ \ref{prop:Rg}.2.\\

When $e$ divides $a$, then $a_0e=a$, so $\xsize{R_q}=\xsize{R}$. So in fact $R_q\cong R$, and $Y^a$ is in the ideal $(q)$ of $\xRo[Y]$. In this case, $q$ is the only member of $\Delta$ annihilating $u$. For if $q_1\in\Delta$ with $q_1(u)=0$ then $(q_1-q)(u)=0$, so $q_1-q\in(q)$. As $\deg(q_1-q)<\deg(q)$, it follows that $q_1=q$. In this case, so if $e\mid a$, we will call the PIR $R$ \textit{\textlabel{neat}{sth:neat}}.\\

When $R$ is not neat, $m\coloneqq a$ mod $e>0$ and $a=(a_0-1)e+m$. In $R_q$, $p\sim\xs{Y}{}^e$, so $\xs{Y}{}^a\sim \xpe{a_0-1}\xs{Y}{}^m$, and hence $R$ $\cong$ $\xRo[Y]/(q,\xpe{a_0-1}Y^m)$. One might regard the rings $R_q$, for the various $q\in\Delta$, as ``neatifications'' of $R$.

\medskip
\begin{prop}\label{prop:neat}Let $R$ be a finite local PIR of type $(a,e,f,p)$ with $a>e$.
\begin{enumerate}[label=\normalfont{(}\normalfont\arabic*)]
\item If $R$ is neat then $R\cong R_g=\xRo[Y]/(g)$, and every uniformizer of $R$ is a zero of a unique polynomial in $\Delta$.
\item When $R$ is not neat, $R\cong\xRo[Y]/(g,\xpe{a_0-1}\,Y^{\xf[.75]a\textup{\xf[1.4]mod\xf[1.4]}e})$.
\item The set $U$ of uniformizers of $R$ equals $\{\xf t(\pi)\pi\mid t\in\xRo[Y]$ and $t(0)\in\xUo\}$.
\item $\xsize{U}=\xsize{\xu}\xsize{k}{}^{a-2}$, if $k$ is the residue field of $R$.
\end{enumerate}
\begin{proof}We just saw (1) and (2). Uniformizers are of the form $\xsum{i=1}{a-1}r_i\pi^i$, with coefficients $r_i\in\rho\subseteq\xRo$ and $r_1\ne0$ (i.e., $r_1\in H_0\subseteq\xUo$). Because $\xsize{\rho}=\xsize{k}$, points (3) and (4) follow.
\end{proof}
\end{prop}

\medskip
Below, we will use the concept of the \textit{resultant} $\res_{\xf Y}(q(Y),w(Y))$ of two polynomials $q(Y)$ and $w(Y)$ in an indeterminate $Y$ over a commutative ring $A$. The resultant is an element of $A$. A good reference for the properties of resultants is \cite{LQ}, Ch.\ III, \S7. We will use their Lemma III.7.3 in particular. Assume that $q(Y)$ is monic of degree $e$. If $B\supseteq A$ is a \textit{\textlabel{splitting ring}{sth:splitr}} for $q$ over $A$, that is, there are $u_i\in B$ such that $q=\xprod{i=1}{e}(Y-u_i)$ in $B[Y]$, then one has $\res_{\xf Y}(q(Y),w(Y))=\xprod{i=1}{e}w(u_i)$ in $B$ (point (4) of Lemma 7.3). Such a $B$ can be constructed as follows. Put $B_1\coloneqq A[X_1]/(q(X_1))$. Then $q(Y)=(Y-\xs{X}_1)q_1(Y)$ for some monic $q_1\in B_1[Y]$ of degree $e-1$. Since $q$ is monic, the natural map $A\to B_1$ is injective. Then put $u_1\coloneqq\xs{X}_1$, adjoin a zero $u_2$ of $q_1$ via $B_2\coloneqq B_1[X_2]/(q_1(X_2))$, and so on. The desired ring $B$ is found in $e-1$ steps.\\

The resultant can be computed directly in $A$ as the value of the \textit{Sylvester determinant} of $q$ and $w$. It is a determinant of order $e+n$, if $w=\xsum{j=0}{n}w_jY^j$. The first row consists of the coefficients of $q$, starting with the leading coefficient, and padded on the right with zeroes. The next $n-1$ rows are repetitions, with each row shifted one position to the right compared to the previous one (and zeroes brought in from the left). Row $n+1$ contains the coefficients $w_n$ thru $w_0$ of $w$, zero-padded. It is repeated $e-1$ times, with shifting. As Lombardi and Quitté note (just prior to Lemma 7.3), because $q$ is monic, this determinant yields the resultant even when $w_n=0$, i.e. when the degree of $w$ is actually $n-1$ (or lower still).\\

As a result, if $d:A\to C$ is a ring homomorphism, $d(\res_{\xf Y}(q(Y),w(Y)))=\res_{\xf Y}(d(q(Y)),d(w(Y)))$ in case $q$ is \textit{monic}, even when $d$ destroys the leading coefficient(s) of $w$. Indeed, $d(q(Y))$ is again monic.\\

Now take $A\coloneqq\xRo[Z]$, where $Z$ is a new indeterminate. Take $q\in\nabla$, so $q\in\xRo[Y]$ is $p\,$-Eisenstein of degree $e$. Let $t\in\xRo[Y]$ with $t(0)\in\xUo$. Put $w(Y)\coloneqq Z-t(Y)Y\in A[Y]$, and denote $\res_{\xf Y}(q(Y),w(Y))$ by $\xh{\lambda(t,q)}$. Factorize $q=\xprod{i=1}{e}(Y-u_i)$ over a splitting ring $B\supseteq A$. Then $\lambda(t,q)=\xprod{i=1}{e}(Z-t(u_i)u_i)$ is monic and of degree $e$ in $Z$. For $j<e$, the coefficient of $Z^j$ is a symmetric and homogeneous polynomial in the $u_i$, so it can be written as a polynomial (without constant term) in the non-leading coefficients of $q(Y)$. The latter are multiples of $p$. Therefore, $\lambda(t,q)\equiv Z^e$ mod $p$. If $d:A\to\xRo$ is given by $Z\mapsto0$, we find $\lambda(t,q)(0)$ $=$ $d(\res_{\xf Y}(q(Y),w(Y)))$ $=$ $\res_{\xf Y}(d(q(Y)),d(w(Y)))$ $=$ $\res_{\xf Y}(q(Y),-t(Y)Y)$ $=$ $\pm\res_{\xf Y}(q(Y),t(Y))\cdot\res_{\xf Y}(q(Y),Y)$ (by point (2) of Lemma 7.3). $\res_{\xf Y}(q(Y),Y)=\pm\,q(0)$ by Lemma 7.3.4, as $Y$ is monic of degree $1$, with $0$ as its zero. Since $q$ is Eisenstein, $\ord_p(q(0))=1$. In the Sylvester matrix for $q$ and $t$, the lower right $e\times e$ minor is lower triangular. As $t(0)\in\xUo$, the only term in the determinant expansion not divisible by $p$ is $1^{\deg(t)}t(0)^e$, the product of the entries on the main diagonal. Hence $\ord_p(\res_{\xf Y}(q(Y),t(Y)))=0$, and $\ord_p(\lambda(t,q)(0))=1$. Therefore, $\lambda(t,q)$ is Eisenstein. Substituting $Y$ for $Z$, this means $\xh{t*q}\coloneqq\lambda(t,q)(Y)$, the \textit{substitution product} of $t$ and $q$, is in $\nabla$ again.\\

This yields (1) of:

\medskip
\begin{prop}\label{prop:resultant}For $q\in\nabla$ and $t\in\xRo[Y]$ with $t(0)\in\xUo$, the following hold.
\begin{enumerate}[label=\normalfont{(}\normalfont\arabic*)]
\item The resultant $\lambda(t,q)$ of $q(Y)$ and $Z-t(Y)Y$ with respect to $Y$ is a polynomial in $Z$ over $\xRo$ that satisfies $t*q\coloneqq\lambda(t,q)(Y)\in\nabla$.
\item The substitution product $t*q$ depends only on the residue class of $t(Y)$ in $\xRo[Y]/(q,\xpe{a_0-1})$.
\item If $q$ is in $\Delta$, so is $t*q$.
\end{enumerate}
\begin{proof}
For (2), $\lambda(t,q)=\lambda(t+c(Y)\cdot q,q)$ for all $c(Y)\in\xRo[Y]$, by Lemma III.7.3.2 in \cite{LQ}. And $\lambda(t+c(Y)\cdot\xpe{a_0-1},q)=\xprod{i=1}{e}(Z-t(u_i)u_i-c(u_i)u_i\xpe{a_0-1})$ if $q=\xprod{i=1}{e}(Y-u_i)$ in $B[Y]$ for a splitting ring $B\supseteq\xRo$. This differs from $\lambda(t,q)=\xprod{i=1}{e}(Z-t(u_i)u_i)$ by $\xpe{a_0-1}$ times a symmetric polynomial in $u_1,\cdots,u_e$ with coefficients in $\xRo[Z]$. The latter can be expressed as a polynomial in the non-leading coefficients of $q(Y)$, hence is a multiple of $p$. But $\xpe{a_0}=0$ in $\xRo$.\\
To see (3), if $q(u)=0$ for a $u\in U$, we have $q=(Y-u)q_1(Y)$ for some $q_1\in R[Y]$. Let $B\supseteq R$ be a splitting ring for $q_1$. Say $q=\xprod{i=1}{e}(Y-u_i)$ in $B[Y]$, with $u_1=u$. Then $t(u)u\in R$ is a zero of $\xprod{i=1}{e}(Z-t(u_i)u_i)=\lambda(t,q)$, and therefore $t*q\in\Delta$.
\end{proof}
\end{prop}

\medskip
The ring $\xh{T}\coloneqq\xRo[Y]/(Y^a)$\index{T@$T$, the ring $\xRo[Y]/(Y^a)$} is local, with maximal ideal $(p,\xs{Y})$. Its units are the $t\in T$ for which $t(0)\in\xUo$. For every $q\in\nabla$, one has $Y^a\in(q,\xpe{a_0-1})\xRo[Y]$ by points (1) and (2) of Prop.\ \ref{prop:neat}, applied to the PIR $R_q/(\xs{Y}{}^a)$. [Or, directly, $q=Y^e-pv$ for some $v\in\xRo[Y]$, hence $Y^{ne}=(q+pv)^n=qw+p^nv^n$, with $w\in\xRo[Y]$, for all $n\in\xN$. Therefore, $Y^{a}=qwY^{a-ne}+p^nv^nY^{a-ne}$ whenever $ne<a\le(n+1)e$. And $a_0-1=n$ for such $a$.] So, due to Prop.\ \ref{prop:resultant}.2, the product $\xh{t*q}$\index{Ttq@$t*q$, the substitution product of a $t\in T_\times\xb$ and a polynomial $q\in\nabla$} is well-defined for $t\in T^\times\xb$. This is, however, generally not a group action. We will define an alternative group structure on the set $T^\times\xb$ to repair this. To spare the eye, below we will simply write $Y$ instead of $\xs{Y}$ in expressions referring to elements of $T$.\\

For $t_1,t_2\in T^\times\xb$, put $\xh{t_1*t_2}\coloneqq t_1(t_2(Y)Y)\cdot t_2(Y)$\index{Ttt@$t_1*t_2$, the substitution product of two elements $t_1,t_2\in T_\times$}, the \textit{substitution product} of $t_1$ and $t_2$. Then $\langle T^\times\xb,*\rangle$ is a group, which we will denote by $\xh{T_\times}\xb$\index{TA@$T_\times$, the set $T^\times\xb$ equipped with the substitution product group structure}. Indeed, $1$ mod $\xpe{a_0}\in \xUo\subseteq T_\times$ is the unit element. We check the remaining requirements. For associativity of $*$, note that $t_1(t_2(t_3(Y)Y)t_3(Y)\cdot Y)\cdot t_2(t_3(Y)Y)t_3(Y)=t_1(t_2(t_3(Y)Y)t_3(Y)Y)\cdot t_2(t_3(Y)Y)\cdot t_3(Y)$ for all $t_1,t_2,t_3\in T_\times\xb$, where the left-hand side evaluates $t_1*(t_2*t_3)$ and the right-hand side $(t_1*t_2)*t_3$.\\

To find the inverse $t_2=\xsum{0\le j<a}{}v_jY^j$ of $t_1=\xsum{0\le i<a}{}u_iY^i$, where $u_i,v_j\in \xRo$, we solve the equation $1=t_1*t_2=(\xsum{0\le i<a}{}u_i\cdot(\xsum{0\le j<a}{}v_jY^{j+1})^i)\cdot(\xsum{0\le j<a}{}v_jY^j)$ recursively, by comparing coefficients of $Y^i\xb$, starting with the constant terms. That gives $v_0=u_0^{-1}$, with the inverse taken in $\xUo$. Let $w_n$ denote the coefficient of $Y^n\xb$ in $\xsum{0\le i<a}{}u_i\cdot(\xsum{0\le j<a}{}v_jY^{j+1})^i$. For $i\ge1$, the coefficient of $Y^i\xb$ on the left-hand side is zero, and on the right it is $\xsum{0\le n\le i}{}\xf w_nv_{i-n}$. Because the $w_n$ for $n\le i$ involve only $v_m$ with $m<i$, which can be assumed to be already known, and since $w_0=u_0$ is a unit, it follows that the coefficients $v_i$ of $t_2$ are uniquely determined.\\

It is clear that this generalizes directly, as follows.

\medskip
\begin{prop}\label{prop:convol}If $A$ is a commutative ring and $a\in\xN$, the substitution product $t_1*t_2=t_1(t_2(Y)Y)\cdot t_2(Y)$ defines a group structure on the set of units of the ring $A[Y]/(Y^a)$.$\hfill\square$
\end{prop}

\medskip
We now verify that $*\xcol T_\times\xb\times\nabla\to\nabla$ is a group action. Take a $q\in\nabla$ and $t_1,t_2\in T_\times\xb$. Let $B\supseteq\xRo$ be a splitting ring for $q$. Say $q=\xprod{i=1}{e}(Y-u_i)$ in $B[Y]$. Then $(t_1*(t_2*q))(Z)=\res_{\xf Y}(\res_{\xf X}(q(X),Y-t_2(X)X),Z-t_1(Y)Y)=\res_{\xf Y}(\xprod{i=1}{e}(Y-t_2(u_i)u_i),Z-t_1(Y)Y)=\xprod{i=1}{e}(Z-t_1(t_2(u_i)u_i)\cdot t_2(u_i)u_i)=\xprod{i=1}{e}(Z-(t_1*t_2)(u_i)\cdot u_i)=\res_{\xf Y}(q(Y),Z-(t_1*t_2)(Y)\cdot Y)=((t_1*t_2)*q)(Z)$.\\

This all but completes the proof of the following.

\medskip
\begin{thm}\label{thm:orbits}Let $R$ be a finite local PIR of signature $(a,e,f,p,\Delta)$, with set of uniformizers $U$ and coefficient ring $\xRo$. Assume $a>e$.
\begin{enumerate}[label=\normalfont{(}\normalfont\arabic*)]
\item * is a group action of\, $T_\times\xb=(\xRo[Y]/(Y^a))^\times\xb$, equipped with the substitution product group structure, on the set $\nabla=\{q\in \xRo[Y]\mid \deg(q)=e$ and $q$ is $p\xf[.5]$-Eisenstein$\}$.
\item $\Delta$ is a union of orbits in $\nabla$ under this action.
\item If $R$ is \ref{sth:neat}, $\Delta$ is the orbit of the Eisenstein polynomial $g$.
\item If $q\in\Delta$ then every $u\in U$ is a zero some some polynomial in the orbit of $q$ under $T_\times\xb$.
\item The polynomials in an orbit all have the same number of zeroes in $R$.
\item The group $T_\times\xb$ equally acts on $U$, through $\xh{t*u}\coloneqq t(u)u$\index{Ttu@$t*u$, the substitution product of a $t\in T_\times\xb$ and a uniformizer $u\in U$}. This action is transitive.
\item If $q\in\Delta$, $u\in U$ and $q(u)=0$, then $(t*q)(t*u)=0$.
\item If $R$ is neat, $U\to\Delta$, $u\mapsto$ the unique $q$ with $q(u)=0$, is an epimorphism of $T_\times\xb$-sets.
\item If $R$ is neat, the isomorphism classes of the finite local PIRs of type $(a,e,f,p)$ are in 1-1 correspondence with the orbits of $T_\times\xb$ in $\nabla$.
\end{enumerate}
\begin{proof}
Point (1) has been proved, (2) follows from Prop.\ \ref{prop:resultant}.3, (7) holds by the proof of that proposition, and (6) by Prop.\ \ref{prop:neat}.3. Checking $(t,u)\mapsto t(u)u$ is a group action is trivial. For (3), if $e$ divides $a$ and $q\in\Delta$, choose a zero $u$ of $q$ in $U$ and write $u=t(\pi)\pi$ for a suitable $t\in T_\times\xb$. By (7), $(t*g)(u)=(t*g)(t*\pi)=0$, because $g(\pi)=0$. Prop.\ \ref{prop:neat}.1 tells us $u$ cannot be a zero of any element of $\nabla$ other than $q$ when $e\mid a$. Hence necessarily $t*g=q$. Item (4) is clear from (7) and (6), item (8) from Prop.\ \ref{prop:neat}.1 and (7). (9) should be obvious from point (3). And $\{u\in U\mid q(u)=0\}\to\{v\in U\mid(t*q)(v)=0\}$, with $u\mapsto t*u$, is well-defined (for $t\in T_\times\xb$) by (7), and bijective by (1) and (7). This gives (5).
\end{proof}
\end{thm}

\medskip
In examples \ref{xmp:2} and \ref{xmp:3}, $\Delta=\{Y^2-2,Y^2+4Y+2\}$ and $\Delta=\{Y^2+2,Y^2+4Y-2\}$, respectively. The remaining Eisensteins are $Y^2\pm2Y+2$ (forming an orbit) and $Y^2\pm2Y-2$ (ditto). The two polynomials in the orbit $\Delta$ each have eight of the 16 uniformizers as roots. As $e\mid a$, there exist four local PIRs of flavour $(6, 2, 1, 2)$ in all. Representatives for the remaining two isomorphism types are $\xZ[Y]/(8,Y^2+2Y+2)$ and $\xZ[Y]/(8,Y^2+2Y-2)$.\\

For example \ref{xmp:4}, $\nabla$ consists of 702 Eisenstein polynomials of degree 2. Half of them have a zero in $R$. These form $\Delta$, the orbit of $g=Y^2+3Y-3$. The remaining half constitute the second orbit. So just two isomorphism classes of type $(3, 2, 3, 3)$ exist. The ring $\xRo[Y]/(Y^2-3Y+3)$ represents the second class.

\begin{xmp}\label{xmp:5}For $a=5$, $e=2$, $f=2$, $p = 2$ and $g(Y)= Y^2+4Y-2$, $\nabla$ contains 192 polynomials. There are 10 orbits, 6 of length 24 and 4 of length 12. The invariant $\Delta$ for the PIR $\xZ[X,Y]/(8,Y^5\xb,g(Y),X^2-X-1)$, of characteristic 8 and order 1024, is the union of the four short orbits. The six size-24 orbits each represent the remaining isomorphism classes of type $(5,2,2,2)$ PIRs.
\end{xmp}

\begin{xmp}\label{xmp:6}For $a=4$, $e=3$, $f=2$, $p = 3$ and $g(Y)= Y^3+6Y^2+3Y+6$, $\Delta=\nabla$. So rings of type $(4,3,2,3)$ are uniquely determined. $\nabla$ consists of 648 polynomials. There are 17 orbits, 8 of length 72 and 9 of length 8. $g$ is in a long orbit, and $Y^3-3$ in a short one. The orbits contained in $\Delta$ are therefore not necessarily all of the same size. There exist 17 isomorphism classes within the neat type $(6\xf[2][=a_0e],3,2,3)$, correlating with the individual orbits.
\end{xmp}

\section{The unramified case}\label{sec:unram}

Here we examine the units and automorphisms of the coefficient ring $\xRo$. The notation is as before. $C_n$ denotes the cyclic group of order $n$, as per usual. We recall from section \ref{sec:pir} that $H_0=\langle z\rangle$ for some unit $z\in H_0\subseteq\rho\subseteq\xRo$. One has $\xUo=H_0\times I_1$, with nested $I_n\coloneqq\{1+rp^n\mid r\in \xRo\}$ for $1\le n<a_0$, similarly to $\xU=H_0\times H_1$. If $\eta\in I_1-I_2$ then $\ord_p(\eta-1)=1$. In the notation of Lemma \ref{lem:elem}, we have $m=1$ and $e=1$, so that the condition $mp>m+e$ is satisfied automatically when $p>2$. In that case, the order of $\eta$ is $p^{\xf[.5]a_0-1}$\xb. Likewise, any $\eta\in I_n-I_{n+1}$ has order $p^{\xf a_0-n}$\xb. Now $I_{a_0-1}\cong k^+\cong(C_p)^f$, for the additive group of $k$ is an elementary abelian $p\xf[.3]$-group. Since $I_n/I_{n+1}\cong k$ for $n<a_0$, it follows that $I_{a_0-2}\cong(C_{p^2})^f$, as the $\eta\in I_{a_0-2}-I_{a_0-1}$ are all of order $p^2$\xb. By induction, $I_1\cong (C_{p^{a_0-1}})^f$.\\

If $p=2$, the $\eta\in I_n-I_{n+1}$ are of order $2^{\xf[.5]a_0-n}$ only when $n>1$. For these $n$, again $I_n\cong(C_{2^{\xf[.5]a_0-n}})^f$. $\xRo$ is of type $(a_0,1,f,2)$. If $a_0=1$, $\xRo$ is just $\xF_{2^f}$, with $I_1=1$. When $a_0=2$, every $\eta\in I_1$ is of the form $1+r\cdot2$ with $r\in\rho$, so of order 1 or 2. So here, $I_1\cong(C_2)^f$. In fact, $v\xcol\eta\mapsto r$ mod $p$ is an isomorphism from $I_1$ to $k^+$. As $k^+$ is elementary abelian, $\xF_2^+=\{0,1\}$ has a complement in $k^+$. When we have $f>1$, there are several complements. No natural choice exists, however. Let $N\le I_1$ be the inverse image under $v$ of a complement of $\xF_2^+$ in $k^+$.\\

For $a_0>2$, let $R'$ be the ring $\xRo/4\xRo$. It is a $(2,1,f,2)$ of the form just discussed, so with group of one-units $I_1'\cong(C_2)^f$. There is a surjection $w\xcol I_1\to I_1'$, given by $1+\xsum{i=1}{a_0-1}r_i\cdot2^i\mapsto1+r_1\cdot2$ (all $r_i\in\rho$). In the notation used above, composing with $I_1'\to I_1'/N\cong\xF_2^+$ gives $I_1\twoheadrightarrow\xF_2^+$. This epimorphism is split by the map $1\mapsto-1=\xsum{i=0}{a_0-1}2^i\in I_1$, since $w(-1)=1+1\cdot2$ corresponds to $1\in\xF_2^+\subseteq k^+$ under the homomorphism $v$, so that $w(-1)\notin N$. Hence $\{\pm1\}$ is a direct summand of $I_1$. For $\eta=1+\xsum{i=1}{a_0-1}r_i2^i\in I_1-I_2$ one has $r_1\ne0$ and $\eta^2\equiv1+(r_1+r_1^2)\cdot4$ mod $8\xRo$. Since $r_1$ is a unit, $\eta^2\in I_3$ iff $\ord_2(1+r_1)>0$ iff $r_1\equiv1$ mod 2 iff $r_1=1$, because $1\in\rho$ and $\rho$ is a system of representatives of the residue classes of $\xRo$ mod 2. Now $J\coloneqq\langle-1\rangle\cdot I_2=\langle-1\rangle\oplus I_2$, because $-1\notin I_2$. So $\xsize{J}=2^{(a_0-2)f+1}$. The elements of $I_2-J$ have their square in $I_2-I_3$, so they are or order $2^{a_0-1}$. This number is therefore the exponent of $I_1$. Hence $I_1\cong(C_{2^{a_0-1}})^n\oplus B$ for some $n\in\xN$ and a subgroup $B$ of exponent $<2^{a_0-1}$. The cyclic group $C_{2^{a_0-1}}$ has $2^{a_0-2}$ generators, hence as many non-generators, so $(C_{2^{a_0-1}})^n$ has $2^{(a_0-1)n}-2^{(a_0-2)n}$ elements of order $2^{a_0-1}$. And $(C_{2^{a_0-1}})^n\oplus B$ has $\xsize{B}$ times as many. Since $I_1$ has $\xsize{I_1-J}$ such elements, it follows that
$\xsize{I_1-J}=2^{(a_0-1)f}-2^{(a_0-2)f+1}=(2^{(a_0-1)n}-2^{(a_0-2)n})\xsize{B}=(2^{(a_0-1)n}-2^{(a_0-2)n})\xf2^{(a_0-1)(f-n)}$. This gives $n=f-1$. Because $I_2=(C_{2^{a_0-2}})^f$ is a subgroup of $I_1$ and $\langle-1\rangle$ a direct summand, we have shown the following.

\medskip
\begin{prop}\label{prop:units}The unit group of $\xRo$ is isomorphic to $C_{p^f-1}\oplus (C_{p^{a_0-1}})^f$ for odd $p$, and to $C_{2^f-1}\oplus(C_{2^{a_0-1}})^{f-1}\oplus C_{2^{a_0-2}}\oplus C_2$ if $p=2$ and $a_0\ge2$. $(\!$When $p=2$ and $a_0=1$, one has $\xRo=k=\xF_{2^f}\xb$, with unit group $C_{2^f-1}$.$)\hfill\square$
\end{prop}

\medskip
$\xRo=\xRb$, with $b$ one of the $f$ zeroes in $\xRo$ of $h$. Cf.\ Prop.\ \ref{prop:omega}.4. The element $b$ mod $p$ generates the residue field $k$ over $\xF_p$. Ring automorphisms of $\xRo$ therefore correspond 1-1 to the ones of $k$.

\medskip
\begin{prop}\label{prop:autos}The automorphism group of $\xRo$ is cyclic of order $f$, generated by the Frobenius automorphism $b\mapsto b\xf'$, where $b\xf'\in\xRo$ is the unique zero of the polynomial $h\in\xZ[X]$ that is congruent to $b^{\xf p}\xb$ mod $p$.
\begin{proof}As $\xsize{\xRo}=\xpef{a_0}=\xsize{(\xZ/\xpe{a_0}\xZ)[X]/(h)}$, the kernel of $(\xZ/\xpe{a_0}\xZ)[X]\twoheadrightarrow\xRo$, with $X\mapsto b$, must be $(h)$. So $\xRo=(\xZ/\xpe{a_0}\xZ)[\xf[.5]b\xf[.5]]\cong(\xZ/\xpe{a_0}\xZ)[X]/(h)$. Because the ring $\xZ/\xpe{a_0}\xZ$ is rigid, the automorphisms of $\xRo$ are given by $b\mapsto b\xf'$, for the $f$ zeroes $b\xf'\in\xRo$ of $h$. Since reduction mod $p$ induces an automorphism of $k=\xF_p(b$ mod $p)=\xF_p(\xs{\beta})$, we have $\Aut(\xRo)\cong\Gal(k/\xF_p)\cong C_f$, and the group is generated by the Frobenius: $b\mapsto$ the $b\xf'\in\{c\in \xRo\mid h(c)=0\}$ for which $b^{\xf p}\xb\equiv b\xf'$ mod $p$.
\end{proof}
\end{prop}

\section{Further notes}\label{sec:further}

A few miscellaneous subjects, such as automorphisms of $R$, quotient rings and subrings, are touched upon here. The phrase \textit{sub}-PIR means a subring that is a PIR. The following notation will be used. $\xh{\nabla/T_\times}\xb$\index{ZNT@$\nabla/T_\times$, the set of orbits of $T_\times\xb$ in $\nabla$} is the set of orbits of $T_\times\xb$ in $\nabla$. For $q\in \nabla$, its orbit is $\xh{T_\times\xb*q}$\index{TOq@$T_\times\xb*q$, the orbit of a polynomial $q\in \nabla$ under $T_\times\xb$}, and $\xh{T_{\times\xb,\xf[.3]q}}\coloneqq\{t\in T_\times\xb\mid t*q=q\}$\index{TSq@$T_{\times\xb,\xf[.3]q}$, the stabilizer in $T_\times\xb$ of a polynomial $q\in \nabla$} denotes the stabilizer. Then $\xsize{T_\times\xb*q}=[T_\times\xb:T_{\times\xb,\xf[.3]q}]$. Likewise, $\xh{T_\times\xb*u}$\index{TOu*@$T_\times\xb*u$, the orbit of a uniformizer $u$ under $T_\times\xb$} is the orbit of a uniformizer $u\in U$ and $\xh{T_{\times\xb,\xf[.3]u}}$\index{TSu@$T_{\times\xb,\xf[.3]u}$, the stabilizer in $T_\times\xb$ of a uniformizer $u$} its stabilizer.\\

Let $G\coloneqq\Aut(R/\xRo)$ be the group of ring automorphisms of $R$ that fix $\xRo$. Since $R=\xRo[\pi]\cong\xRo[Y]/(Y^a,g)$, a $\sigma\in G$ is determined by a zero $\sigma(\pi)$ of $g$ in $R$. If $u\in R$ and $g(u)=0$, then $u\in U$, so $u^a=0$. The map $\xRo[Y]/(Y^a,g)\to R$ with $\pi=\xs{Y}\mapsto u$ is therefore well-defined. It is surjective because $u$ is a uniformizer. Hence, by finiteness of $R$, it is an automorphism. We obtain a homomorphism $T_{\times\xb,\xf[.3]g}\to G$, given by $t\mapsto(\pi\mapsto t*\pi)$. For by Th.\ \ref{thm:orbits}.7, $t*\pi=t(\pi)\pi$ is a zero of $t*g=g$. For every $u\in U$ there is a $t\in T_\times\xb$ with $u=t*\pi$, so if $g(u)=0$ then $\pi$ is a zero of $t^{-1}*g$ (with the inverse of $t$ taken in $T_\times\xb$, of course). In case $e\mid a$, $g$ is the unique Eisenstein annihilating $\pi$, so $t^{-1}*g=g$, and thus $t\in T_{\times\xb,\xf[.3]g}$. So for \ref{sth:neat} $R$, we have $T_{\times\xb,\xf[.3]g}\twoheadrightarrow G$, with kernel $T_{\times\xb,\xf[.3]\pi}$. This gives:

\medskip
\begin{prop}\label{prop:relautos}The relative ring automorphisms of $R$ over $\xRo$ correspond to the zeroes of $g$ in $R$. For neat $R$, the group $\Aut(R/\xRo)$ is isomorphic to $T_{\times\xb,\xf[.3]g}/T_{\times\xb,\xf[.3]\pi}$ and it is of order $\xsize{\{u\in R\mid g(u)=0\}}=\xsize{U}/\xsize{\Delta}$.
\begin{proof}
According to Th.\ \ref{thm:orbits}.6, $U=T_\times\xb*\pi$. Hence $\xsize{U}=[T_\times\xb:T_{\times\xb,\xf[.3]\pi}]$. If $e\mid a$, by Th.\ \ref{thm:orbits}.3 we have $\xsize{\Delta}=[T_\times\xb:T_{\times\xb,\xf[.3]g}]$, so $|\xb\Aut(R/\xRo)\xf[.3]|=\xsize{T_{\times\xb,\xf[.3]g}/T_{\times\xb,\xf[.3]\pi}}=\xsize{U}/\xsize{\Delta}$.
\end{proof}
\end{prop}

\medskip
Regarding finite modules, we have:

\begin{prop}\label{prop:modules}Every finitely generated $R$-module $M$ is a finite direct sum of modules of the form $R/\xp^i$ with $0<i\le a$.
\begin{proof}
Let $D$ be a DVR that maps onto $R$ (Th.\ \ref{thm:main}.1). Then $D$ is a principal ideal domain, and $M$ is a finitely generated $D$-module. As is well-known, $M$ is then a finite direct sum of cyclic $D$-modules.
\end{proof}
\end{prop}

\medskip
Let $0<i<a$, and let $R'$ be the quotient ring $R/\xp^i$. The residue field of $R'$ equals $k$, but one has $p^{\xceil{i/e}}\xb=0$ in $R'$. So the coefficient ring $\xRo'$ of $R'$ may be a proper quotient of $\xRo$. If $i\le e$ then $p=\xo{\pi}^{\xf[.3]e}=0$ in $R'$, so $R'$ is equicharacteristic, hence isomorphic to $k[Y]/(Y^i)$. When $i>e$, $R'$ has signature $(i,e,f,p,\Delta')$, say, with $\Delta\subseteq\Delta'\subseteq\nabla$. Indeed, the Eisensteins having roots in $R$ also have roots in $R'$. It may occur that $\Delta'=\Delta$. [A case in point is the $R$ from ex.\ \ref{xmp:4}, seen as a quotient of its $R_g=\xRo[Y]/(Y^2+3Y-3)$.] The natural map $T\twoheadrightarrow T'\coloneqq\xRo'[Y]/(Y^i)$ induces a surjection $T^\times\xb\twoheadrightarrow T'{}^\times\xb$, which is a homomorphism $T_\times\xb\to T'_\times\xb$ wrt.\ the substitution product. The orbits of the group $T'_\times\xb$ in $\nabla'$ therefore coincide with those of the group $T_\times\xb$.

\medskip
\begin{prop}\label{prop:quotients}Let $e<i<a$, and let $R'=R/\xp^i$ be the quotient ring.
\begin{enumerate}[label=\normalfont{(}\normalfont\arabic*)]
\item $R'$ is of signature $(i,e,f,p,\Delta')$, with $\Delta\subseteq\Delta'\subseteq\nabla$, and it has coefficient ring $\xRo'=\xRo/(p^{\xceil{i/e}})\xb$.
\item $\mathrm{nat}\xcol T_\times\xb\twoheadrightarrow T'_\times\xb$ is a group homomorphism, where $T'=\xRo'[Y]/(Y^i)$.
\item $\nabla/T'_\times\xb=\nabla/T_\times\xb$. That is, $T'_\times\xb$ and $T_\times\xb$ have the same orbits in $\nabla$.$\hfill\square$
\end{enumerate}
\end{prop}

\medskip
Subrings of $R$ are local rings, and they contain $R_0=\xZ/\xpe{a_0}\xZ$. By Prop.\ \ref{prop:basics}.8, the residue field of a subring is always a subfield of $k$. Not every subring is a PIR. E.g., $R_0+\xp^i$ has maximal ideal $(p)+\xp^i$. If $i\ge e$, $p$ is an element of minimal $\ord_\xp$, and otherwise $\pi^i$ is. So when $i>e$ and $e\nmid i$, or $i<e$ and $i\nmid e$, this ring cannot be a PIR.\\

If $R'$ is a sub-PIR of $\xRo$ and $k'\subseteq k$ its residue field, $R'$ is generated over $R_0$ by the elements that are equal to their $\xsize{k'}$-th power. In view of the uniqueness of Hensel lifting, $R'$ is determined as a subring of $\xRo$ by the subfield $k'$ of $k$, which in turn is determined by the divisor $f'=[k':\xF_p]$ of $f=[k:\xF_p]$.\\

Sub-PIRs of $R$ have a sub-PIR of $\xRo$ as coefficient ring. Confining ourselves to sub-PIRs $R'$ that contain $\xRo$, we are dealing with rings of some type $(a',e',f,p)$. Uniformizers $\pi'$ of $R'$ are $\sim\pi^d$ in $R$ for some $d\in\xN$. Then $(\pi')^{e'}\sim p\sim\pi^e$, hence $e=e'd$ and $a'=\xceil{a/d\xf[1.5]}$. If $\Delta'\subseteq\nabla'$ denotes the Delta invariant of $R'$, again nat$\xcol T\to T'=\xRo[Y]/(Y^{a'})$ induces an epimorphism $T_\times\xb\twoheadrightarrow T'_\times\xb$, and the $T'_\times\xb$-orbits in $\nabla'$ agree with the $T_\times\xb$-orbits. Let $\Delta'_R\subseteq\nabla'$ be the set of Eisensteins of degree $e'$ having a root in $R$. It is a union of orbits, a number of which make up $\Delta'$. In case $e'\mid a'$, $\Delta'$ simply is one of these orbits, and the remaining orbits in $\Delta'_R$ belong to other sub-PIRs of type $(a',e',f,p)$.\\

Recall the surjection $\psi\xcol\xO_L\to R$ constructed in section \ref{sec:structure}. For $R'$, there is of course a similar map. But we can obtain one from $\psi$ by choosing a $\gamma'\in\xO_L$ with $\psi(\gamma')=\pi'$. Let $L'\coloneqq K(\gamma')$. Then clearly $\psi\xb[2.5]\upharpoonright\xb\xO_{L'}\xcol\xO_{L'}\twoheadrightarrow R'$. The sub-PIRs of $R$ containing $\xRo$ therefore match the subextensions $K\subseteq L'$ of $K\subseteq L$.\\

A general sub-PIR $R'$ of $R$ is an $(a',e',f',p)$, again with $e'\mid e$. We have $f'\mid f$ as well, as the residue field of $R'$ is a subfield of the one of $R$. Here one has $\xO_{L'}\twoheadrightarrow R'$ for a subextension $\xQ\subseteq L'$ of $\xQ\subseteq L$. For the finite local PIR $R'$, the role of $K$ is played by the field $K\cap L'$.

\medskip
\begin{prop}\label{prop:subrings}In relation to the field $L$ and the homomorphism $\psi\xcol\xO_L\twoheadrightarrow R$ defined in section \ref{sec:structure}, the sub-PIRs of $R$ are precisely the images under $\psi$ of the rings of integers $\xO_{L'}$ of the intermediate fields $\xQ\subseteq L'\subseteq L$.$\hfill\square$
\end{prop}

\medskip
If $A$ is a commutative ring, $\xp$ a prime ideal with $\xp\ne\xp^2\xb$, and $e\in\xN$, let $\nabla=\nabla_{\xb[1.2]A,\xp,e}$ be the set of the $q(Y)\in A[Y]$ of degree $e$ that are Eisenstein wrt.\ $\xp$. That is, $q$ is monic, all coefficients except the leading one are in $\xp$, and the constant term is in $\xp-\xp^2\xb$. Such $q$ are irreducible by the usual argument: if $q=q_1q_2$ with $q_1,q_2\in A[Y]$ monic and not constant, then $\xo{q}_i=Y^{\deg(q_i)}$ in $(A/\xp)[Y]$, so $q_1(0),q_2(0)\in\xp$, contradicting $q(0)\notin\xp^2\xb$. $A[Y]_\times\xb\coloneqq\{\xf t\in A[Y]\mid t(0)\in A^\times\}$ is a monoid under the substitution product. It operates on $\nabla$ through $t*q\coloneqq\res_{\xf Y}(q(Y),Z-t(Y)Y)(Y)$. The arguments of section \ref{sec:delta}, showing that $t*q$ is again Eisenstein clearly remain valid in the present context. This renders:

\medskip
\begin{prop}\label{prop:monoid}With the above notation, $*\xcol A[Y]_\times\xb\times \nabla\to \nabla$ is a monoid action on the set of $\xf[.5]\xp$-Eisenstein polynomials of degree $e$ over $A$.$\hfill\square$
\end{prop}

\medskip
We do not know whether the substitution multiplication, an ad-hoc term, and groups like $T_\times\xb$ have already been described in the literature. In the simple case $a=e=2$, one has $\xRo=k$, and $T_\times\xb$ is $k^+\xb[2]\rtimes\xu$, the semidirect product of the normal subgroup $\{x\xs{Y}+1\mid x\in k\}$ and the subgroup $\xu$. The size and complexity of $T_\times\xb$ grow rapidly. For instance, for $(a,e,f,p)=(3,2,2,2)$ one finds $T_\times\xb\cong(C_4\times C_4\times(((C_4\times C_2\times C_2) \rtimes C_2)\rtimes C_2))\rtimes C_3$, of order $3\cdot4^5=3072$. In general, $\xsize{T_\times}\xb=\xsize{\xu}\xsize{k}{}^{a_0a-1}$. By Prop.\ \ref{prop:resultant}.2, though, $T_\times\xb$ may be replaced by the smaller group $((\xRo/\xpe{a_0-1})[Y]/(Y^a))_\times\xb$ of the same kind, if $a>e$.

\phantomsection
\addcontentsline{toc}{section}{Appendix on orders in number fields}\section*{Appendix on orders in number fields}\label{sec:appx}

Certain properties of orders in algebraic number fields that were used in the text are covered here. The material is more or less well-known. I came across Lemma \ref{appl:invert} among my (disorganized) notes, without a reference to its origin.\\

Let $K$ be an algebraic number field, i.e.\ a finite extension of $\xQ$, and $\xO_K$ the ring of integers. An \textit{order} in $K$ is a subring $A\subseteq K$ that is a free of rank $n=[K:\xQ]$ over $\xZ$. Equivalently, $A$ is a subring of finite index of $\xO_K$ (\cite{ST}, Th.\ 2.2). Every order is a Noetherian ring of Krull dimension 1.\\

A \textit{fractional ideal} of $A$ is an $A$-submodule $I$ of $K$ such that $aI\subseteq A$ for some $0\ne a\in A$. $aI$ is then a (non-zero) ideal of $A$, and $I=a^{-1}\cdot aI$ is a finitely generated $A$-module. As $K$ is the quotient field of $A$ and $\xO_K$ is finitely generated as an $A$-module (indeed, as a $\xZ$-module), $\xO_K$ is a fractional ideal of $A$.\\

For fractional ideals $I$ and $J$, the sum $I+J$, product $IJ$, and \textit{ideal quotient} $(I:J)=\{x\in K\mid xJ\subseteq I\}$ are again fractional ideals. $I$ is called \textit{invertible} if $IJ=A$ for some fractional ideal $J$ of $A$. In that case, one has $J=(A:I)$, and we denote this (unique) inverse $J$ by $I^{-1}$. Non-zero ideals of $A$ are called invertible if they are invertible as a fractional ideals.\\

The \textit{conductor} of the order $A$ is $(A:\xO_K)=\{x\in K\mid x\xO_K\subseteq A\}$. It is contained in $A$, for $1\in\xO_K$. It is easy to see that the conductor is the largest ideal of $\xO_K$ that is contained in $A$.

\medskip
\begin{appl}\label{appl:lem1}If $I,J,H\subseteq K$ are fractional ideals of $A$ and $I$ is invertible, one has $I(J\cap H)=IJ\cap IH$.
\begin{proof}
$\subseteq$ is clear. In particular, $I^{-1}(IJ\cap IH)\subseteq I^{-1}IJ\cap I^{-1}IH$. As $I^{-1}I=A=II^{-1}$, multiplication with $I$ gives $IJ\cap IH\subseteq I(J\cap H)$.
\end{proof}
\end{appl}

\medskip
\begin{appl}\label{appl:invert}Let $I\subseteq K$ be a fractional ideal of $A$. The following are equivalent.
\begin{enumerate}[label=\normalfont{(}\normalfont\arabic*)]
\item $I$ is invertible.
\item For every ideal $0<J\le A$ there exists an $x\in K$ such that $A=xI+J$.
\end{enumerate}
\begin{proof}
(1) $\Rightarrow$ (2). We may assume that $J<A$. Let $\xp_1,\cdots,\xp_m$ be the different prime ideals of $A$ that contain $J$. They correspond to the minimal primes of the Noetherian ring $A/J$, so there are finitely many of them. If $1\le i\le m$, the product $\xprod{1\le j\le m,\xf j\ne i}{}\xf\xp_j$ of the ideals $\xp_j$ for $j\ne i$ cannot be $\subseteq\xp_i$. Therefore, $I^{-1}\cdot\xprod{1\le j\le m,\xf j\ne i}{}\xf\xp_j\nsubseteq I^{-1}\cdot\xprod{1\le j\le m}{}\xf\xp_j$. Taking $x_i$ in $I^{-1}\xprod{1\le j\le m,\xf j\ne i}{}\xf\xp_j$ but not in $I^{-1}\xprod{1\le j\le m}{}\xf\xp_j$, each $x_i$ is in $I^{-1}$, hence so is $x\coloneqq\xsum{1\le i\le m}{}x_i$. I.e.\ $xI\subseteq A$. Suppose $xI+J<A$. Then $xI+J\subseteq\xp_i$ for some $i$, say $i=1$. For $j>1$ we have $Ix_j\subseteq II^{-1}\xprod{1\le l\le m,\xf l\ne j}{}\xf\xp_l\subseteq\xp_1$, so $x_j\in I^{-1}\xp_1$. Since $xI\subseteq xI+J\subseteq\xp_1$, $x\in I^{-1}\xp_1$ too. So $x_1$ is both in $I^{-1}\xp_1$ and $I^{-1}\xprod{2\le j\le m}{}\xf\xp_j$. As $\xp_1$ and $\xprod{2\le j\le m}{}\xf\xp_j$ are comaximal in $A$, their intersection is equal to their product. By Lemma \ref{appl:lem1}, $x_1\in I^{-1}\xprod{1\le j\le m}{}\xf\xp_j$, clashing with the choice of $x_1$. So $xI+J=A$.\\
(2) $\Rightarrow$ (1). Clearly, $J=(A:I)I$ is a non-zero ideal of $A$. Take an $x\in K$ with $A=xI+J$. Then $x\in(A:I)$, so $xI\subseteq J$. Hence $A=(A:I)I$.
\end{proof}
\end{appl}

\medskip
If $J\ne0$ is an ideal of $A$ and $0\ne a\in J$, $a$ is integral over $\xZ$, so $f(a)=0$ for some $f\in\xZ[X]$ with non-zero constant term $m=f(0)$. Then $m\in J$. As $A$ is free of rank $n$ over $\xZ$, $A/mA$ is finite of size $m^n$. A fortiori, the index $[A:J]$ is finite. This number is known as the \textit{absolute norm} or \textit{ideal norm} of $J$.

\medskip
\begin{appp}\label{appp:absnorm}If $I\subseteq K$ is an invertible fractional ideal of $A$ and $0<J\le A$ an ideal, then $A/J\cong I/IJ$ as $A$-modules. For every pair of non-zero ideals $J,H$ of $A$ with at least one invertible, $[A:JH]=[A:J][A:H]$ . In particular, if $J$ is an invertible ideal of $A$ and $m\in\xN$, one has $[A:J^m]=[A:J]^m$.
\begin{proof}
We may assume $J<A$. By Lemma \ref{appl:invert}, there is an $x\in K^\times\xb$ such that $xI+J=A$. We have $A/J=(xI+J)/J\cong xI/(xI\cap J)$. Because $xI$ and $J$ are comaximal in $A$, $xI\cap J=xIJ$. And $x\cdot\xb[1.5]\xcol I\to xI$ is an $A$-isomorphism sending $IJ$ onto $xIJ$. For the second statement, by the first, $[A:H]=[J:JH]$ (if $J$ is invertible) $=[A:JH]/[A:J]$. The third claim follows at once.
\end{proof}
\end{appp}

\medskip
\begin{appl}\label{appl:lem2}If $\xp$ is a maximal ideal of the order $A$ and there is only one prime ideal, $\xP$, of $\xO_K$ that lies over $\xp$, then $(\xO_K)_\xp$ is a discrete valuation ring.
\begin{proof}
$(\xO_K)_\xp=S^{-1}\xO_K$, where $S\coloneqq A-\xp$, is a localization of a Dedekind domain, hence Dedekind. Its maximal ideals are the $\mathfrak{M}\cdot(\xO_K)_\xp$ with $\mathfrak{M}\in\max(\xO_K)$ such that $\mathfrak{M}\xf[.5]\cap\xf[.5]S=\varnothing$, i.e.\ $\mathfrak{M}\xf[.5]\cap\xf[.5]A\subseteq\xp$, i.e.\ $\mathfrak{M}=\xP$. So $(\xO_K)_\xp$ is a local Dedekind domain (with maximal ideal $\xP(\xO_K)_\xp$), and therefore a DVR.
\end{proof}
\end{appl}

\medskip
\begin{appp}\label{appp:prinv}Let $\xp$ be a maximal ideal of the order $A\subseteq\xO_K$ and $\xP$ an $\xO_K$-prime lying over $\xp$. Then the following conditions are equivalent.
\begin{enumerate}[label=\normalfont{(}\normalfont\arabic*)]
\item $\xp$ is an invertible ideal.
\item For every maximal ideal $\mathfrak{q}$ of $A$, the $A_\mathfrak{q}$-module $\xp A_\mathfrak{q}$ is monogenous. 
\item The maximal ideal $\xp A_\xp$ of $A_\xp$ is principal.
\item $A_\xp$ is a discrete valuation ring.
\item $A_\xp=(\xO_K)_{\xP}$, as subrings of $K$.
\item $\xp+(A:\xO_K)=A$, that is, $\xp$ is prime to the conductor.
\item $\xP=\xp\xO_K$ and $A/\xp\cong\xO_K/\xP$.
\end{enumerate}
\begin{proof}
(1) $\Leftrightarrow$ (2) follows from (a) $\Leftrightarrow$ (d) of Th.\ 4 in \cite{B}, \S II.5.6, which applies to any commutative ring $A$ and an $A$-submodule $\xp$ of $S^{-1}A$ with $\xp S^{-1}A=S^{-1}A$, if $S$ is a multiplicative subset of $A$ that doesn't contain any zero-divisors. $K$ is the quotient field of $A$, and we use $S\coloneqq A-\{0\}$. Clause (d) also specifies $\xp$ is finitely generated, but that is redundant here since $A$ is Noetherian.\\
(2) $\Leftrightarrow$ (3) is trivial, for if $\mathfrak{q}\ne\xp$ is a maximal ideal of $A$ then $\xp A_\mathfrak{q}=A_\mathfrak{q}$ is monogenous over $A_\mathfrak{q}$.\\
(3) $\Leftrightarrow$ (4) holds by \cite{B}, \S VI.3.6, (a) $\Leftrightarrow$ (d) of Prop.\ 9.\\
(4) $\Leftrightarrow$ (5) since $A_\xp\subseteq(\xO_K)_{\xP}$ as subrings of $K$, $(\xO_K)_{\xP}$ is a DVR, and the DVRs of $K$ are the $\subseteq$-maximal proper subrings of $K$.\\
Assuming (5), every $x\in\xO_K\subseteq(\xO_K)_{\xP}$ is of the form $a/s$ with $a\in A$ and $s\in A-\xp$. As $\xO_K$ is f.g.\ over $A$, there is an $s\in A-\xp$ with $s\xO_K\subseteq A$. I.e.\ $s\in(A:\xO_K)$, and (6) follows because $\xp$ is a maximal ideal of $A$. Conversely, given (6), $1=s+a$ for some $s\in(A:\xO_K)$ and $a\in\xp$. Then $s\in A-\xp$, and for every $x\in\xO_K$ and $t\in\xO_K-\xP$ one has $x/t=(sx)/(st)\in A_\xp$. Hence $(\xO_K)_{\xP}\subseteq A_\xp$, i.e.\ $A_\xp=(\xO_K)_{\xP}$. As $x=sx+ax\equiv sx$ mod $\xp\xO_K$, the map $A\to\xO_K/\xp\xO_K$ is surjective. The kernel contains $\xp$ and clearly cannot be larger. So $\xO_K/\xp\xO_K\cong A/\xp$, a field, and therefore $\xp\xO_K=\xP$. Thus (6) implies (7).\\
We now show (7) $\Rightarrow$ (4). Let $\mathfrak{f}=(A_\xp:(\xO_K)_\xp)$ be the conductor. By Lemma \ref{appl:lem2}, $(\xO_K)_\xp$ is a DVR, so $\mathfrak{f}=\xP^{m}\cdot(\xO_K)_\xp$ for some $m\ge0$. Because $\xsize{A/\xp}=\xsize{\xO_K/\xP}$, the natural map $A/\xp\rightarrowtail\xO_K/\xP$ is a bijection, and we have $\xO_K=A+\xP=A+\xp\xO_K=A+\xp(A+\xp\xO_K)=A+\xp^2\xO_K$. By induction, $\xO_K=A+\xp^m\xO_K$, so every $y\in(\xO_K)_\xp$ is of the form $(a+\xsum{i=1}{l}b_ix_i)/s$ for some $a\in A$, $l\in\xN$, $b_i\in\xp^m$, $x_i\in\xO_K$ and $s\in A-\xp$. But the $b_i$ are in $\mathfrak{f}$, hence $y=a/s+\xsum{i=1}{l}b_i\cdot(x_i/s)\in A_\xp$. So $(\xO_K)_\xp=A_\xp$ (and $m=0$ after all). Thus $A_\xp$ is a DVR.
\end{proof}
\end{appp}

\medskip
We include one of various results that go by the name of Gauss' Lemma. Recall that a polynomial $g=\xsum{i=0}{m}a_iX^i$ over a ring $A$ is \textit{primitive} if its \textit{content ideal} $\xsum{i=0}{m}a_iA$ is equal to $A$ (or, if you like, if $g(A)=A$).

\medskip
\begin{appp}[Gauss' Lemma]\label{appp:GL}
For commutative rings $A\subseteq B$, the following are equivalent.
\begin{enumerate}[label=\normalfont{(}\normalfont\arabic*)]
\item $A$ is integrally closed in $B$.
\item If $f\in A[X]$ factorizes as $f=gh$ in $B[X]$ with $g$ and $h$ monic, then $g$ and $h$ are in $A[X]$.
\end{enumerate}
Under these conditions one has:
\begin{enumerate}[resume*]
\item If $f\in A[X]$ factorizes as $f=gh$ in $B[X]$ with $g\in A[X]$ primitive, then $h$ is in $A[X]$.
\end{enumerate}
If $B$ is a domain, properties $(1)$ and $(2)$ are equivalent to:
\begin{enumerate}[resume*]
\item Every monic, irreducible $f\in A[X]$ remains irreducible in $B[X]$.
\end{enumerate}
\begin{proof}
(1) $\Rightarrow$ (2): take an extension $C\supseteq B$ that is a \ref{sth:splitr} for both $g$ and $h$. Then we have $g = \xprod{i=1}{n}(X-x_i)$ and $h = \xprod{j=1}{m}(X-y_j)$ in $C[X]$ for suitable $x_i,y_j\in C$. As $f$ is monic and $f(x_i)=0=f(y_j)$ for all $i$ and $j$, the $x_i$ and $y_j$ are integral over $A$. The coefficients of $g$ and $h$ are polynomials in the $x_i$ and $y_j$, respectively, so they are integral over $A$ too. As these coefficients are in $B$, by (1) they are in $A$.\\
To see (2) $\Rightarrow$ (1), if $f(b)=0$ for a $b\in B$ and some monic $f\in A[X]$, then $f = (X-b)g$ for a polynomial $g\in B[X]$. Then $g$ is monic, so by (2) $X-b$ and $g$ are in $A[X]$. In particular, $b\in A$.\\
For (1) $\Rightarrow$ (3), let $g=\xsum{i=0}{m}a_iX^i\in A[X]$, with $\xsum{i=0}{m}a_ic_i=1$ for certain $c_i\in A$, and $h=\xsum{j=0}{n}b_jX^j\in B[X]$. By a theorem of Kronecker's (\cite{LQ}, Th.\ III.3.3), also known as the Dedekind Prague theorem, for each $0\le j\le n$ the products $a_ib_j$ ($0\le i\le m$) are integral over the subring of $A$ generated (over $\xZ$) by the coefficients of $gh=f$, hence over $A$ itself. But then the $b_j=\xsum{i=0}{m}(a_ib_j)c_i$ are integral over $A$, so $h\in A[X]$ due to (1).\\
Now let $B$ be a domain. (2) $\Rightarrow$ (4): if $f = gh$ in $B[X]$ and $g$ and $h$ are not units, their leading coefficients are inverse units in $B$. So we may assume $g$ and $h$ are monic. By (2), they are in $A[X]$ and neither is a unit. Contradiction. For (4) $\Rightarrow$ (1): if $f(b) = 0$ with $b\in B$ and $f\in A[X]$ monic, factorize $f$ as $f_1\cdots f_r$ with the $f_i$ monic and irreducible in $A[X]$. It is clear that this is possible. Since $B$ is a domain, $f_i(b) = 0$ for some $i$. By (4), $f_i$ is irreducible in $B[X]$. But $X-b$ divides $f_i$ in $B[X]$. Hence $f_i = X-b$. As $f_i\in A[X]$, we have $b\in A$.
\end{proof}
\end{appp}

\medskip

The Kummer-Dedekind theorem (c.f.\ Th.\ 8.2 in \cite{ST}) remains valid when one replaces the base ring $\xZ$ by the ring of integers of a number field $K$. To be precise, let $g\in\xO_K[X]$ be monic and irreducible. By (1) $\Rightarrow$ (4) of Prop.\ \ref{appp:GL}, $g$ is irreducible in $K[X]$. Let $L=K[X]/(g)=K(\gamma)$, where $\gamma\coloneqq\xs{X}\in L$ is a zero of $g$. Then $\gamma\in\xO_L$, the ring of integers of $L$. Let $B$ be the order $\xO_K[\gamma]$, and $\xp$ a maximal ideal of $\xO_K$. We have $\xo{g}=\xprod{i=1}{r}\xo{g}_i^{\xf e_i}$ in $(\xO_K/\xp)[X]$ for some $r,e_1,\cdots,e_r\in\xN$ and suitable monic $g_i\in\xO_K[X]$ for which the $\xo{g}_i$ are different and irreducible in $(\xO_K/\xp)[X]$. Put $\xP_i=\xp B+g_i(\gamma)B$ for $1\le i\le r$. For each $i$, write $g=g_iq_i+h_i$ in $\xO_K[X]$, with remainder $h_i\in\xp\cdot\xO_K[X]$ of degree less than $\deg(g_i)$. The ``quotient'' polynomials $q_i$ are then monic too.

\medskip
\begin{appt}[Kummer-Dedekind]\label{appt:KD}With notation as above:
\begin{enumerate}[label=\normalfont{(}\normalfont\arabic*)]
\item The $\xP_i$ are the prime ideals of $B$ that lie over $\xp$.
\item $\xprod{i=1}{r}\xP_i^{\xf e_i}\subseteq\xp B$.
\item Equality holds in $(2)$ iff every $\xP_i$ is invertible.
\item A particular $\xP_i$ is invertible iff $e_i=1$ or $h_i\in\xp\xf[.3]\xO_K[X]-\xp^2\xf[.3]\xO_K[X]$.
\end{enumerate}
\begin{proof}
We follow Stevenhagen's proof of Th.\ 8.2 in \cite{ST}, which represents the case $K=\xQ$. $B/\xP_i\cong\xO_K[X]/((g)+(g_i)+\xp[X])\cong(\xO_K/\xp)[X]/(\xo{g}_i)$ is a field, so $\xP_i\in\max(B)$. If $\xp\subseteq\xP\in\max(B)$, in $B/\xP$ one has $0=g(\xo{\gamma})=\xprod{i=1}{r}g_i(\xo{\gamma})^{\xf e_i}$, hence $g_i(\gamma)\in\xP$ for some $i$. Therefore, $\xP=\xP_i$, and (1) results.\\
Every element of $\xprod{i=1}{r}\xP_i^{\xf e_i}$ is congruent mod $\xp B$ to a multiple of $\xprod{i=1}{r}g_i(\gamma)^{e_i}\equiv g(\gamma)$ mod $\xp B$, hence $\equiv0$, which gives (2).\\
As $\xO_K$ is a Dedekind domain, the ideal $\xp$ is invertible, so $\xp\xf[.3]J=\xO_K$ for some fractional $\xO_K$-ideal $J\subseteq K$. So if $\xprod{i=1}{r}\xP_i^{\xf e_i}=\xp B$, we have $(\xprod{i=1}{r}\xP_i^{\xf e_i})\cdot JB=B$, and by definition the $\xP_i$ are invertible ideals of $B$. Conversely, if each $\xP_i$ is invertible, by Prop.\ \ref{appp:absnorm} $[B:\xprod{i=1}{r}\xP_i^{\xf e_i}]=\xprod{i=1}{r}[B:\xP_i]^{\xf[.3]e_i}$. Now $[B:\xP_i]=\xsize{B/\xP_i}=\xsize{\xO_K/\xp}^{\xf[.3]\deg(g_i)}\xb$, so $[B:\xprod{i=1}{r}\xP_i^{\xf e_i}]=\xsize{\xO_K/\xp}{}^{\xf[.3]\deg(g)}\xb$. But $\xsize{B\xf[.3]/\xp B}=\xsize{(\xO_K/\xp)[X]/(\xo{g})}$ also equals $\xsize{\xO_K/\xp}{}^{\xf[.3]\deg(g)}\xb$. This proves (3).\\
Before proceeding with the proof of (4), we note that, given some $1\le i\le r$ and $t\in(\xO_K)_{\xp}[X]$, one has:\\
\-\qquad\quad$t(\gamma)\in\xP_iB_{\xP_i}\;\Leftrightarrow\;\;\xo{g}_i\mid\xo{t}$ \;in $(\xO_K)_{\xp}/\xp\xf[.3](\xO_K)_{\xp}[X]\cong(\xO_K/\xp)[X]$\hfill($*$)\\
Indeed, under the composition $B_{\xP_i}\to B_{\xP_i}/\xP_iB_{\xP_i}\xrightarrow{\sim}(\xO_K)_{\xp}/\xp\xf[.3](\xO_K)_{\xp}[X]/(\xo{g}_i)$, the element $t(\gamma)$ maps to the congruence class of $\xo{t}\in(\xO_K)_{\xp}/\xp\xf[.3](\xO_K)_{\xp}[X]$ mod $\xo{g}_i$.\\
Now for (4), first assume that $e_i=1$. Then $\xo{g}_i\nmid\xo{q}_i$ in $(\xO_K/\xp)[X]$, so by ($*$) we have $q_i(\gamma)\notin\xP_i$. So $q_i(\gamma)\in B_{\xP_i}^{\xf[.3]\times}\xb$. But $0=g(\gamma)= g_i(\gamma)q_i(\gamma)+h_i(\gamma)$ and $h_i\in\xp\xf[.3]\xO_K[X]$, hence $g_i(\gamma)\in\xp B_{\xP_i}$. As $(\xO_K)_{\xp}$ is a DVR, $\xp\xf[.3](\xO_K)_{\xp}$ is a principal ideal. It follows that $\xP_iB_{\xP_i}=\xp B_{\xP_i}=\xp\xf[.3](\xO_K)_{\xp}\cdot B_{\xP_i}$ is principal too. By (3) $\Rightarrow$ (1) of Prop.\ \ref{appp:prinv}, $\xP_i$ is invertible. Now assume $h_i\notin\xp^2\xf[.3]\xO_K[X]$. Let $\pi\in\xp$ be a uniformizer, i.e.\ $\xp\xf[.3](\xO_K)_{\xp}=\pi(\xO_K)_{\xp}$. Then $h_i=\pi h_0$ in $(\xO_K)_{\xp}[X]$ for some $h_0\in(\xO_K)_{\xp}[X]-\xp\xf[.3](\xO_K)_{\xp}[X]$. So $0\ne\xo{h}_0\in(\xO_K)_{\xp}/\xp\xf[.3](\xO_K)_{\xp}[X]\cong(\xO_K/\xp)[X]$. Because $\deg(h_0)<\deg(g_i)$, $\xo{g}_i\nmid\xo{h}_0$ in $((\xO_K)_{\xp}/\xp\xf[.3](\xO_K)_{\xp})[X]$ and hence $h_0(\gamma)\in B_{\xP_i}^{\xf[.3]\times}\xb$ by ($*$). But then we have $\pi=h_0(\gamma)^{-1}h_i(\gamma)=-h_0(\gamma)^{-1}g_i(\gamma)q_i(\gamma)\in g_i(\gamma)B_{\xP_i}$, hence $\xP_iB_{\xP_i}=(g_i(\gamma))$ is again principal. For the converse, assume that $e_i\ge2$ and $h_i\in\xp^2\xf[.3]\xO_K[X]$. Then $\xo{g}_i\mid\xo{q}_i$, hence $q_i(\gamma)\in\xP_i$. And, with $\pi$ and $h_0$ as above, $h_0\in\xp\xf[.3](\xO_K)_{\xp}[X]=\pi(\xO_K)_{\xp}[X]$, so $h_0(\gamma)=\pi x$ in $B_{\xP_i}$ for some $x\in B_{\xP_i}$. Let $y\coloneqq q_i(\gamma)/\pi\in L$. Then $y g_i(\gamma)=g_i(\gamma)q_i(\gamma)/\pi=-h_i(\gamma)/\pi=-h_0(\gamma)=-\pi x\in\xP_iB_{\xP_i}$. Hence $y\xP_iB_{\xP_i}=y\pi B_{\xP_i}+yg_i(\gamma)B_{\xP_i}=q_i(\gamma)B_{\xP_i}+\pi x B_{\xP_i}\subseteq\xP_iB_{\xP_i}$. Therefore, $\xP_iB_{\xP_i}$ is a faithful $B_{\xP_i}\xb[.75][\xf[.5]y\xf[.5]]$-module which is finitely generated over $B_{\xP_i}$. Ergo, $y$ is integral over $B_{\xP_i}$ (see e.g. \cite{B}, \S V.1.1, Th.\ 1). If $y\in B_{\xP_i}$\xb, we have $q_i(\gamma)v(\gamma)=\pi u(\gamma)$ in $B$ for some $u,v\in\xO_K[X]$ with $v(\gamma)\notin\xP_i$. So $g$, being the minimal polynomial of $\gamma$ over $K$, divides $q_iv-\pi u$ in $K[X]$, hence also in $\xO_K[X]$ because $g$ is monic. Thus, in $(\xO_K/\xp)[X]$, $\xo{q}_i\xo{g}_i=\xo{g}\mid\xo{q}_i\xo{v}$, so that $\xo{g}_i\mid\xo{v}$. But this contradicts $v(\gamma)\notin\xP_i$. Therefore, $y\notin B_{\xP_i}$. So $B_{\xP_i}$ is not integrally closed and it cannot be a DVR. By Prop.\ \ref{appp:prinv}, the ideal $\xP_i$ is not invertible.
\end{proof}
\end{appt}

\medskip
If $\xP_i$ is invertible, by Prop.\ \ref{appp:prinv} $\tilde{\xP_i}\coloneqq\xP_i\xO_L$ is the unique prime of $\xO_L$ that lies over $\xP_i$. The $B$-prime $\xP_i$ is then called \textit{regular}, and otherwise \textit{singular}. The singular $\xP_i$ may split and/or ramify on extension to $\xO_L$, or lead to a further extension of the residue field. They contain the conductor $(B:\xO_L)$, so there are only finitely many of them in $\max(B)$, for the combined $\xp\in\max(\xO_K)$. When all $\xP_i$ are invertible, $\xp\xO_L=\xprod{i=1}{r}\tilde{\xP_i}{}^{\xb e_i}$ is the prime ideal factorization of the extension $\xp\xO_L$ of $\xp$ to $\xO_L$. In this case, the order $B$ is called \textit{regular above} $\xp$.\\

\phantomsection

\phantomsection
\addcontentsline{toc}{chapter}{Index of symbols}
\printindex


\addcontentsline{toc}{chapter}{References}

\begin{thebibliography}{17}

\bibitem{B} N. Bourbaki, Commutative Algebra, Chapters 1-7, Springer Verlag, New York $\cdot$ Berlin $\cdot$ Heidelberg, 2nd printing, 1989.

\bibitem{CR} Coefficient rings, The Stacks project, \href{https://stacks.math.columbia.edu/tag/0326}{Definition 0326}.

\bibitem{CO} I.S. Cohen, On the Structure and Ideal Theory of Complete Local Rings, \textit{Trans. AMS} \underline{59}-1 (1946), 54-106.

\bibitem{HL} Hensel's Lemma, Wikipedia article, \href{https://en.wikipedia.org/wiki/Hensel's_lemma}{Hensel's Lemma}.

\bibitem{HU} T.W. Hungerford, On the Structure of Principal Ideal Rings, \textit{Pacific J. of Math.} \underline{25} (1968), 543-547.

\bibitem{LQ} H. Lombardi and C. Quitté, Commutative Algebra: Constructive Methods, (Algebra and Applications; 20), Springer Verlag, New York $\cdot$ Heidelberg $\cdot$ Berlin, 2015.

\bibitem{SE} J.-P. Serre, Local Fields, (Graduate Texts in Mathematics; 67), Springer Verlag, New York $\cdot$ Heidelberg $\cdot$ Berlin, 1979.

\bibitem{ST} P. Stevenhagen, The arithmetic of number rings. In: Algorithmic Number Theory: Lattices, Number Fields, Curves and Cryptography, Cambridge University Press, Cambridge $\cdot$ New York $\cdot$ Oakleigh, (2008), 209-266. (The article is available from the \href{https://www.math.leidenuniv.nl/~psh/ANTproc/08psh.pdf}{Leiden University}.)

\end{thebibliography}
\end{document}